\documentclass[twoside]{article}
\usepackage[preprint]{aistats2027}

\usepackage[utf8]{inputenc}
\usepackage{natbib}
\usepackage[colorlinks=true,linkcolor=blue!50!black,citecolor=blue!50!black,urlcolor=blue!50!black]{hyperref}
\usepackage{url}
\usepackage{booktabs}
\usepackage{amsfonts}
\usepackage{nicefrac}
\usepackage[expansion=false]{microtype}
\usepackage{xcolor}
\usepackage{amssymb,amsthm,mathtools}
\usepackage{graphicx}
\usepackage{algorithm}
\usepackage[noend]{algorithmic}
\usepackage{enumitem}
\usepackage{tikz}
\usetikzlibrary{shadows, positioning, calc, decorations.pathreplacing}
\usetikzlibrary{arrows.meta}
\usepackage[capitalize,noabbrev]{cleveref}

\newtheorem{theorem}{Theorem}
\newtheorem{lemma}[theorem]{Lemma}
\newtheorem{proposition}[theorem]{Proposition}
\newtheorem{corollary}[theorem]{Corollary}
\theoremstyle{definition}
\newtheorem{definition}[theorem]{Definition}
\newtheorem{remark}[theorem]{Remark}

\crefname{theorem}{Theorem}{Theorems}
\crefname{lemma}{Lemma}{Lemmas}
\crefname{proposition}{Proposition}{Propositions}
\crefname{corollary}{Corollary}{Corollaries}
\crefname{definition}{Definition}{Definitions}
\crefname{remark}{Remark}{Remarks}
\Crefname{theorem}{Theorem}{Theorems}
\Crefname{lemma}{Lemma}{Lemmas}
\Crefname{proposition}{Proposition}{Propositions}
\Crefname{corollary}{Corollary}{Corollaries}
\Crefname{definition}{Definition}{Definitions}
\Crefname{remark}{Remark}{Remarks}

\newcommand{\defeq}{\triangleq}
\newcommand{\R}{\mathbb{R}}
\newcommand{\E}{\mathbb{E}}
\newcommand{\eps}{\varepsilon}
\newcommand{\Ht}{\widetilde{H}}
\newcommand{\gt}{\widetilde{g}}
\newcommand{\fstar}{f^{\star}}
\newcommand{\Ftilde}{\widetilde{F}}
\newcommand{\taubar}{\bar{\tau}}
\newcommand{\tmax}{\tau_{\max}}
\newcommand{\grad}{\nabla f}
\newcommand{\hess}{\nabla^2 f}
\newcommand{\lmin}{\lambda_{\min}}
\newcommand{\norm}[1]{\left\lVert #1 \right\rVert}
\newcommand{\inner}[2]{\left\langle #1, #2 \right\rangle}
\newcommand{\Ov}{O}                 % Overlap count
\newcommand{\Cov}{C_{\mathrm{ov}}}

\newcommand{\HBFGS}{B}              % BFGS curvature approx

\newcommand{\cmark}{\ensuremath{\checkmark}}
\newcommand{\xmark}{\ensuremath{\text{--}}}

\begin{document}

\runningtitle{Second-Order Optimization with Asynchronous Curvature}
\runningauthor{Chayti, Jaggi}

\twocolumn[
\aistatstitle{Second-Order Optimization with Asynchronous Curvature}
\aistatsauthor{El Mahdi Chayti \And Martin Jaggi}
\aistatsaddress{Machine Learning and Optimization Laboratory (MLO), EPFL}
]

\begin{abstract}
Second-order methods make fast progress per iteration, but forming and factorizing the Hessian often dominates their wall-clock time. We study \emph{Split-Client}, a task-parallel design in which a curvature client factorizes Hessians in the background while a gradient client takes cubic-regularized steps with the most recently published matrix, so that curvature cost appears as a delay rather than as blocking time. For delayed exact Hessians, we prove a first-order stationarity bound that depends on the average age of the curvature rather than its maximum, using an overlap count and a weighted Hardy inequality that allows any bounded-error initial matrix; the analysis extends to subsampled and uniformly inexact curvature. Under a fixed-cost model, the resulting wall-clock bound matches the scaling of optimally tuned Lazy Hessian without choosing a reuse period. In experiments on two LIBSVM datasets, an implementation decreases the objective earlier in wall-clock time than Lazy and Vanilla cubic Newton, with the largest gains at intermediate accuracy.
\end{abstract}

% ============================================================================
\section{Introduction}
\label{sec:intro}
% ============================================================================

Newton's method and cubic regularization use curvature to make fast progress per iteration and are robust to ill-conditioning, with fast local convergence under suitable local assumptions \citep{conn2000trust, nesterov2006cubic, cartis2011adaptive}. These properties make them attractive for machine learning, where training objectives are often nonconvex and poorly conditioned.

Their wall-clock cost limits their use. For a finite sum of $n$ terms in dimension $d$, a gradient costs $\mathcal{O}(nd)$, forming the Hessian costs $\mathcal{O}(nd^2)$, and \emph{factorizing} it, for instance by Cholesky or eigendecomposition to solve the cubic subproblem, costs $\mathcal{O}(d^3)$. One curvature job therefore costs $q_{\mathrm{grad}} = \Theta(d + d^2/n)$ gradient evaluations in this arithmetic model; whether assembly or factorization dominates depends on $d/n$ and on the implementation. Standard implementations are \emph{synchronous}: every iteration waits for the factorization to complete. We call this the synchronization barrier.

\paragraph{Regime of applicability.}
We target problems where (i) the full $d \times d$ Hessian fits in memory, roughly $d \in [10^2, 10^4]$ on commodity hardware, and (ii) $d$ is large enough that a curvature job is slow compared with a gradient pass. This regime covers scientific computing, small-to-moderate statistical models, polynomial and kernel features, graphical models, and the inner loops of bilevel and AutoML procedures. The range is illustrative: depending on memory, hardware, sparsity and structure, larger problems may be better served by Hessian-free methods (Newton-CG via Hessian-vector products) or structured approximations (K-FAC \citep{martens2015optimizing}, Shampoo \citep{gupta2018shampoo}); the Split-Client abstraction could be applied there with a different cost model for the curvature step, which we leave to future work.

% --- Figure 1 -----------------------------------------------------------------
\begin{figure}[t]
\centering
\resizebox{\columnwidth}{!}{\begin{tikzpicture}[
    font=\sffamily\scriptsize,
    >=Latex,
    node distance=0.4cm,
    box/.style={draw, rounded corners=1pt, minimum height=4.5mm, inner sep=0pt, align=center,
        drop shadow={opacity=0.15, shadow xshift=0.5pt, shadow yshift=-0.5pt}},
    G/.style={box, fill=green!15, draw=green!40!black, text=green!30!black},
    H/.style={box, fill=blue!12, draw=blue!40!black, text=blue!40!black},
    D/.style={box, fill=red!12, draw=red!50!black, text=red!50!black},
    C/.style={box, fill=gray!10, draw=gray!60!black},
    pics/clock/.style={code={
        \draw[thick, fill=white] (0,0) circle (1.8mm);
        \draw[very thin] (0,0) -- (90:1.3mm);
        \draw[very thin] (0,0) -- (0:0.9mm);
        \foreach \x in {0,90,180,270} \draw (\x:1.5mm) -- (\x:1.8mm);
        \node[yshift=2.6mm, scale=0.4] at (0,0) {\textbf{!}};}},
    lab/.style={font=\sffamily\footnotesize\bfseries, inner sep=0pt, anchor=west}]
\def\wG{0.6cm}\def\wS{0.8cm}\def\wH{1.4cm}\def\wD{1.6cm}
% Panel A
\begin{scope}[yshift=1.6cm]
  \node[lab] at (-4, 0.5) {(A) Vanilla Cubic};
  \node[anchor=west, scale=0.8, text=gray] at (-4, 0.1) {Blocking: Total Time = $T_{grad} + T_{hess} + T_{dec}$};
  \node[G, minimum width=\wG] (aG) at (-3.7,-0.5) {G};
  \node[H, right=0.3cm of aG, minimum width=\wH] (aH) {Hessian};
  \node[D, right=0.3cm of aH, minimum width=\wD] (aD) {Decomp};
  \node[C, right=0.3cm of aD, minimum width=\wS] (aS) {Step};
  \draw[->] (aG) -- (aH); \draw[->] (aH) -- (aD); \draw[->] (aD) -- (aS);
  \path ($(aG.east)!0.5!(aH.west)$) pic {clock};
  \path ($(aH.east)!0.5!(aD.west)$) pic {clock};
  \path ($(aD.east)!0.5!(aS.west)$) pic {clock};
\end{scope}
% Panel B
\begin{scope}[yshift=0cm]
  \node[lab] at (-4, 0.5) {(B) Lazy Hessian};
  \node[anchor=west, scale=0.8, text=gray] at (-4, 0.1) {Blocking on Refresh, Fast on Reuse};
  \node[G, minimum width=\wG] (bG1) at (-3.7,-0.5) {G};
  \node[H, right=0.3cm of bG1, minimum width=\wH] (bH) {Hessian};
  \node[D, right=0.3cm of bH, minimum width=\wD] (bD) {Decomp};
  \node[C, right=0.3cm of bD, minimum width=\wS] (bS1) {Step};
  \node[G, right=0.4cm of bS1, minimum width=\wG] (bG2) {G};
  \node[C, right=0.3cm of bG2, minimum width=\wS, label={[scale=0.6, text=gray]below:Reuse $\Ht$}] (bS2) {Step};
  \draw[->] (bG1) -- (bH); \draw[->] (bH) -- (bD); \draw[->] (bD) -- (bS1);
  \draw[->] (bS1) -- (bG2); \draw[->] (bG2) -- (bS2);
  \path ($(bG1.east)!0.5!(bH.west)$) pic {clock};
  \path ($(bH.east)!0.5!(bD.west)$) pic {clock};
\end{scope}
% Panel C
\begin{scope}[yshift=-2cm]
  \node[lab] at (-4, 0.8) {(C) Split-Client (Ours)};
  \node[anchor=west, scale=0.8, text=gray] at (-4, 0.4) {Non-blocking: Curvature hidden by gradients};
  \node[anchor=east, font=\sffamily\tiny\bfseries, text=gray] at (-4.2, -0.2) {Gradient Worker};
  \draw[dotted, gray!50] (-4.2, -0.6) -- (5.5, -0.6);
  \node[G, minimum width=\wG] (cG1) at (-3.7,-0.2) {G};
  \node[C, right=0.15cm of cG1, minimum width=\wS] (cS1) {Step};
  \node[G, right=0.15cm of cS1, minimum width=\wG] (cG2) {G};
  \node[C, right=0.15cm of cG2, minimum width=\wS] (cS2) {Step};
  \node[G, right=0.15cm of cS2, minimum width=\wG] (cG3) {G};
  \node[C, right=0.15cm of cG3, minimum width=\wS] (cS3) {Step};
  \draw[->] (cG1) -- (cS1); \draw[->] (cS1) -- (cG2);
  \draw[->] (cG2) -- (cS2); \draw[->] (cS2) -- (cG3); \draw[->] (cG3) -- (cS3);
  \node[anchor=east, font=\sffamily\tiny\bfseries, text=gray] at (-4.2, -1.0) {Curvature Worker};
  \node[H, minimum width=\wH] (cH) at (-3.3,-1.0) {Hessian};
  \node[D, right=0.15cm of cH, minimum width=\wD] (cD) {Decomp};
  \draw[->] (cH) -- (cD);
  \draw[->, thick, dashed, red!70!black] (cD.east) -- ++(0.2,0) |- node[pos=0.75, below, font=\tiny, text=red] {New $\Ht$ arrives} ($(cS3.south)+(0,0.05)$);
  \draw[decorate, decoration={brace, amplitude=3pt, mirror}, gray]
    ($(cH.south west) + (0,-0.1)$) -- ($(cD.south east) + (0,-0.1)$)
    node[midway, below=3pt, font=\tiny, text=black] {Total Curvature Delay $\tau$};
\end{scope}
\begin{scope}[yshift=-4cm]
  \draw[->, thick, gray] (-3.7, 0) -- (5.5, 0) node[right] {Wall-Clock Time};
  \node[box, fill=green!15, draw=green!40!black, minimum height=3mm, minimum width=4mm] at (-2, -0.6) {};
  \node[right, scale=0.8] at (-1.8, -0.6) {Gradient};
  \node[box, fill=blue!12, draw=blue!40!black, minimum height=3mm, minimum width=4mm] at (0, -0.6) {};
  \node[right, scale=0.8] at (0.2, -0.6) {Hessian Eval};
  \node[box, fill=red!12, draw=red!50!black, minimum height=3mm, minimum width=4mm] at (2.2, -0.6) {};
  \node[right, scale=0.8] at (2.4, -0.6) {Factorization};
  \path (4, -0.6) pic {clock};
  \node[right, scale=0.8] at (4.2, -0.6) {Idle/Block};
\end{scope}
\end{tikzpicture}}
\caption{Wall-clock comparison of cubic regularization strategies. \textbf{(A) Vanilla:} The gradient worker waits for Hessian computation and factorization at every step. \textbf{(B) Lazy:} Reuses each factorization for $m$ steps, amortizing its cost, but still blocks at every refresh. \textbf{(C) Split-Client (Ours):} Runs curvature in a background process. The gradient worker starts immediately with a cheap initial matrix $H_0$ and proceeds continuously using the latest published $\Ht$, so the curvature cost (Hessian \emph{and} factorization) appears as a delay rather than as blocking time.}
\label{fig:compare-strategies-improved}
\end{figure}

\paragraph{Beyond the arithmetic cost model.}
Most analyses of efficient second-order methods measure \emph{arithmetic} complexity. Lazy Hessian \citep{doikov2023lazy} amortizes curvature computation and factorization by reusing each factorization for $m$ steps. The reuse period is a parameter chosen before the run: in the fixed-cost model of \cref{sec:wall_clock_comparison}, the period that minimizes the wall-clock upper bound is close to the duration of a curvature job measured in iterations, which depends on hardware, linear-algebra implementation and problem size. The analysis does not treat curvature latency as a \emph{delay} that governs convergence.

\paragraph{The Split-Client framework.}
We study a \textbf{task-parallel} alternative. A \textbf{gradient client} evaluates gradients and takes cubic steps continuously, while a \textbf{curvature client} computes and factorizes Hessians in the background. The gradient client never waits: at iteration $k$ it uses the most recently published matrix, whose age $\tau_k$ is determined by the running system rather than chosen in advance. Two difficulties arise. The ages vary over time, so bounds in terms of the maximum age $\tmax$ are conservative. And until the first curvature job completes, the method must use an initial matrix that may be far from the Hessian.

\paragraph{Contributions.}
\begin{enumerate}[leftmargin=*,itemsep=0.1em]
    \item \textbf{Two-timeline model.} A model separating gradient iterations from wall-clock time, in which a curvature job (Hessian assembly and factorization) has its own duration and the age of the curvature in use is a system variable (\cref{sec:framework}).
    \item \textbf{Average-age analysis.} For cubic regularization with exact delayed Hessians and a regularization proportional to the observed age, a first-order bound that depends on the average age $\taubar$ rather than on $\tmax$ (\cref{thm:exact}). It rests on an \emph{Overlap Count}, which bounds how often a step is reused when job durations grow at a bounded rate, and on a \emph{weighted Hardy inequality}, which controls the startup phase: any symmetric initial matrix with bounded error may be used, at the cost of an additive term uniform in the startup duration. Under a common fixed-cost model with job duration $D$, the wall-clock upper bound scales as $\sqrt{D}$, the exponent of optimally tuned Lazy Hessian, without a reuse period chosen in advance; against Lazy with the arithmetic choice $m = d$, the ratio of the compared leading upper-bound terms is $\tfrac{1}{\sqrt2}(\sqrt{d/D} + \sqrt{D/d})$, up to method-specific constants (\cref{cor:lazy_md}).
    \item \textbf{Inexact curvature.} Extensions to curvature and gradient errors proportional to recent step sizes, with causal tolerances realized by adaptive subsampling (\cref{thm:unified}), and to a uniformly bounded curvature error through an inflated schedule (\cref{thm:robust}).
    \item \textbf{Experiments.} On two LIBSVM datasets, an implementation of Split-Client decreases the objective earlier in wall-clock time than Lazy and Vanilla cubic Newton, with the largest gains at intermediate accuracy (\cref{sec:experiments}).
\end{enumerate}

% ============================================================================
\section{Related Work}
\label{sec:related}
% ============================================================================

\paragraph{Efficient second-order methods.}
Inexact Newton \citep{dembo1982inexact, eisenstat1996inexact} and cubic regularization \citep{nesterov2006cubic, cartis2011adaptive} provide global guarantees, and subsampled cubic regularization \citep{xu2020newton, kohler2017sub} reduces per-iteration arithmetic; all remain synchronous. \citet{chayti2023unified, chayti2022optauxinfo} study inexact curvature for iteration complexity. \citet{kamzolov2023cubic} analyze cubic-regularized Quasi-Newton methods under a uniform error bound, with a quadratic inexactness correction, for convex and star-convex objectives; the uniform-error result of \cref{thm:robust} concerns the nonconvex asynchronous setting.

\paragraph{Lazy and reused curvature.}
Reusing curvature goes back to the Shamanskii method \citep{shamanskii1967modification, fan2013shamanskii}; \citet{doikov2023lazy} reuse it for several steps. Lazy Hessian amortizes Hessian computation and factorization through reuse in an arithmetic-cost model, with a reuse period chosen in advance. Split-Client instead models curvature latency as an asynchronous delay determined by the running system; \cref{sec:wall_clock_comparison} compares the two.

\paragraph{Asynchronous and distributed optimization.}
Asynchrony is well studied for \emph{first-order} methods: lock-free analyses \citep{recht2011hogwild, lian2015async}, robustness to arbitrary delays \citep{cohen2021asynchronous, mishchenko2022asynchronous}, delay-adaptive step sizes \citep{koloskova2022sharper, wu2022delay}, and optimal time complexities \citep{tyurin2023optimal}. There the delayed quantity is a gradient and the setting is data-parallel. Here the delayed quantity is a Hessian factorization, whose $\mathcal{O}(d^3)$ cost is charged separately from the $\mathcal{O}(nd)$ gradient, and the parallelism is between two tasks on the same objective. Distributed second-order methods \citep{wang2018giant, crane2019dingo, islamov2021newton, islamov2022distributed, safaryan2022fednl} split the data across devices and target communication, which is complementary.

\begin{table*}[t]
\centering\footnotesize
\setlength{\tabcolsep}{4pt} % Reduces the default horizontal padding between columns
\caption{Positioning relative to related cubic-regularized methods. \emph{Delay model}: curvature latency is modeled as an asynchronous delay. \emph{Avg-delay}: the bound depends on the average age $\taubar$ rather than $\tmax$. \emph{No reuse period}: none must be chosen in advance.}
\label{tab:positioning}
\begin{tabular}{lcccc}
\toprule
Method & Non-block. & Delay model & Avg-delay & No reuse period \\
\midrule
Vanilla CR \citep{nesterov2006cubic} & \xmark & \xmark & n/a & \cmark \\
Inexact CR \citep{chayti2023unified} & \xmark & \xmark & n/a & \cmark \\
Lazy Hessian \citep{doikov2023lazy} & \xmark & \xmark & n/a & \xmark \\
QN-CR \citep{kamzolov2023cubic}      & \xmark & \xmark & n/a & \cmark \\
\textbf{Split-Client (ours)}         & \cmark & \cmark & \cmark & \cmark \\
\bottomrule
\end{tabular}
\end{table*}

% ============================================================================
\section{The Split-Client Framework}
\label{sec:framework}
% ============================================================================

We minimize a twice continuously differentiable $f : \R^d \to \R$, possibly nonconvex, with $d$ moderate so the full Hessian fits in memory.

\begin{definition}[Two-timeline model]
\label{def:timeline}
Time is measured in two ways: by \emph{gradient iterations} $k = 0, 1, 2, \ldots$, at each of which the gradient client evaluates a gradient $g_k$ at $x_k$, computes a step $s_k$ and sets $x_{k+1} = x_k + s_k$; and by \emph{wall-clock time}, normalized so that one full iteration (gradient evaluation and cubic solve) takes one unit. Curvature jobs are processed by a FIFO-serial queue on integer boundaries: $a_0 = 0$, job $i$ runs on $[a_i, b_i]$ with $a_{i+1} = b_i$ and duration $\Delta_i \defeq b_i - a_i \in \{1, 2, \ldots\}$, reads $x_{a_i}$, and publishes $\hess(x_{a_i}) + E_i$, which is used from iteration $b_i$ on. At iteration $k$ with $a_i \le k < b_i$, the in-use curvature is $\Ht_k = \hess(x_{a_{i-1}}) + E_{i-1}$ for $i \ge 1$, and an initial symmetric matrix $H_0$ for $i=0$. The \emph{delay} is $\tau_k \defeq k - a_{i-1}$ for $i \ge 1$. Physical completion times may be rounded into this protocol, in which case (A3) below is imposed on the resulting integer durations.
\end{definition}

\paragraph{The unit-iteration convention.}
The normalization in \cref{def:timeline}, one full iteration per unit of wall-clock time, is a modeling convention: writing $g_k$ for the physical duration of iteration $k$, including its cubic solve with the published factorization, we take $g_k = 1$. It makes comparisons controlled, since two methods are compared as parallel worlds performing the same per-iteration work on the same objective and differing only in how curvature is obtained. Arithmetic analyses such as Lazy Hessian already include factorization in the cost of a curvature update; our model differs in treating curvature computation as a task that overlaps with gradient steps, whose latency appears as the age of the matrix in use. \cref{app:phys_time} gives an elapsed-time bound valid for arbitrary $g_k$.

\begin{figure*}[t]
\centering
\begin{tikzpicture}[
    font=\sffamily\footnotesize, >=Latex, x=0.75cm, y=0.6cm
]
% --- Wall-clock axis ---
\draw[->, thick] (-0.3, 0) -- (15, 0) node[right] {wall-clock $k$};
% Tick marks at integers (gradient iterations)
\foreach \k in {0,1,...,14} {
  \draw (\k, -0.05) -- (\k, 0.05);
}
% Iteration labels (sparse)
\node[below=1pt, scale=0.85] at (0,0) {$0$};
\node[below=1pt, scale=0.85] at (3,0) {$a_1$};
\node[below=1pt, scale=0.85] at (6,0) {$a_2$};
\node[below=1pt, scale=0.85] at (10,0) {$a_3$};
\node[below=1pt, scale=0.85] at (14,0) {$a_4$};

% --- Curvature jobs as colored intervals (lane below) ---
\definecolor{job1}{RGB}{198,219,239}
\definecolor{job2}{RGB}{158,202,225}
\definecolor{job3}{RGB}{107,174,214}
\definecolor{job4}{RGB}{66,146,198}

\draw[fill=job1, draw=blue!50!black] (0, -1.6) rectangle (3, -1) node[midway, scale=0.85] {Job 0};
\draw[fill=job2, draw=blue!50!black] (3, -1.6) rectangle (6, -1) node[midway, scale=0.85] {Job 1};
\draw[fill=job3, draw=blue!50!black] (6, -1.6) rectangle (10, -1) node[midway, scale=0.85] {Job 2};
\draw[fill=job4, draw=blue!50!black] (10, -1.6) rectangle (14, -1) node[midway, scale=0.85] {Job 3};

\node[anchor=east, scale=0.85] at (-0.3, -1.3) {curvature jobs};

% Job durations brace
\draw[decorate, decoration={brace, amplitude=2pt, mirror}, gray!70]
  (3, -1.7) -- (6, -1.7) node[midway, below=2pt, scale=0.8, text=gray] {$\Delta_1$};
\draw[decorate, decoration={brace, amplitude=2pt, mirror}, gray!70]
  (6, -1.7) -- (10, -1.7) node[midway, below=2pt, scale=0.8, text=gray] {$\Delta_2$};

% --- Highlight an iteration k inside Job 2, and its delay ---
\def\kpos{8}
\draw[red!70!black, very thick] (\kpos, 0) -- (\kpos, 0.5);
\node[red!70!black, scale=0.9, above=1pt] at (\kpos, 0.5) {$k$};

% Curvature in use at k: from a_1 (start of Job 1, the previous completed job)
\draw[red!70!black, very thick, dashed] (3, 0) -- (3, 0.5);
\node[red!70!black, scale=0.85, above=1pt] at (3, 0.5) {$a_{i-1}$};

% Delay arrow
\draw[<->, thick, red!70!black] (3, 0.7) -- (\kpos, 0.7) node[midway, above, scale=0.85] {$\tau_k = k - a_{i-1}$};

% --- Overlap window: future k' that include k in their delay window ---
% Job 2 = current job; future k' in jobs 2 (after k) and 3 see k in their delay window.
\draw[decorate, decoration={brace, amplitude=2pt}, orange!70!black, thick]
  (\kpos+0.05, 0.4) -- (14, 0.4) node[midway, above=2pt, scale=0.85, text=orange!70!black]
  {overlap window: $O_k \le \Delta_i + \Delta_{i+1} - 1$};

\end{tikzpicture}
\caption{The two-timeline model. Curvature jobs (blue blocks) run in serial on the wall-clock axis; gradient iterations are integer ticks above. At iteration $k$ inside job $i = 2$, the in-use curvature is from job $i-1 = 1$ (started at $a_1$), giving delay $\tau_k = k - a_1$. The Overlap Lemma (\cref{lem:overlap}) bounds the count of future iterations $k' > k$ that still include $k$ in their delay window: such $k'$ must lie in jobs $\le i+1 = 3$, so $O_k \le \Delta_i + \Delta_{i+1} - 1$.}
\label{fig:timeline}
\end{figure*}

The delay $\tau_k$ is a system variable determined by the relative cost of gradient and curvature computation, not a hyperparameter. At each iteration, the gradient client computes the cubic-regularized step using the fresh gradient $g_k$ (the exact $\grad(x_k)$, or a sampled estimate in \cref{sec:structured}) and the latest published curvature:
\begin{equation}
    s_k \in \arg\min_{s \in \R^d} \big\{ \inner{g_k}{s} + \tfrac{1}{2}\inner{\Ht_k s}{s} + \tfrac{\rho_k}{6}\norm{s}^3 \big\},
    \label{eq:cubic_step}
\end{equation}
then $x_{k+1} = x_k + s_k$, with $\rho_k > 0$ a regularization parameter (full pseudocode in \cref{app:algorithm}). Every curvature matrix, $H_0$ included, is symmetric; replacing an approximation $B$ by $(B + B^\top)/2$ never increases its operator-norm distance to a symmetric Hessian and imposes no positivity condition. The minimizer is global, and deterministic rate statements take expectations over the output index, conditional on the run.

\paragraph{Cost model.} With $D$ the duration of one curvature job (Hessian and factorization) in units of one full iteration: \textbf{Vanilla} blocks every step (per-step $1+D$); \textbf{Lazy} ($m$-reuse) blocks every $m$ steps (per-step $1+D/m$); \textbf{Split-Client} never blocks (per-step $1$, with stale curvature); until the first job completes it uses the initial matrix $H_0$. Split-Client stores two $d \times d$ matrices (double buffering), about $64$\,MB in double precision at $d=2000$.

% ============================================================================
\section{Theory: Exact Delayed Curvature}
\label{sec:theory_exact}
% ============================================================================

We first treat $E_i = 0$ for all $i \ge 0$: the curvature client returns exact Hessians, only delayed.

\subsection{Assumptions}

\begin{enumerate}[label=(A\arabic*),itemsep=0em]
    \item \textbf{Smoothness:} $\hess$ is $L$-Lipschitz for some $L > 0$: $\norm{\hess(x) - \hess(y)} \le L\norm{x-y}$. (If the Hessian is constant, any positive upper bound may be used.)
    \item \textbf{Bounded below:} $f(x) \ge \fstar$; let $F_0 \defeq f(x_0) - \fstar$.
    \item \textbf{Bounded growth of job durations:} there exists a finite $\gamma \ge 1$ such that $\Delta_{i+1} \le \gamma\Delta_i$ for every $i$.
    \item \textbf{Bounded delay:} $\tau_k \le \tmax < \infty$.
    \item \textbf{Initialization:} the initial matrix $H_0$ is symmetric and satisfies $\norm{H_0 - \hess(x_0)} \le \delta_0$; it may be indefinite, for instance a multiple of the identity or zero.
\end{enumerate}

Assumption (A3) bounds increases in successive job durations; decreases need not be bounded. In the exact-curvature analysis, $\gamma$ enters only through the overlap constant $\Cov = \gamma(1+\gamma)$ of \cref{lem:overlap}; the structured analysis of \cref{sec:theory_inexact} also uses $\Delta_i \le \gamma^2\Delta_{i-2}$, from two applications of (A3). Dense factorization at fixed $d$ has predictable runtime (\cref{app:cost_stability}).

\begin{remark}[On the cost-stability constant]
\label{rem:cov}
We do \emph{not} assume $\gamma$ is close to $1$. All results hold for any finite $\gamma$, but our sufficient choice of the multiplier $c$ of the regularization schedule introduced in \cref{sec:overlap} grows linearly in $\Cov$ (\cref{lem:work}), so it requires a valid bound on $\gamma$, and the constants of the rates carry this dependence. Taking $\gamma \to 1$ gives $\Cov \to 2$ and the smallest constants.
\end{remark}

\subsection{The Overlap Count}
\label{sec:overlap}

\Cref{def:timeline} defines the delay $\tau_k$ only once the first curvature job has completed. Write $h \defeq a_2$, and extend the delay to every iteration through the dispatch index of the matrix in use, which we call the \emph{age}:
\begin{equation}
    d_k \defeq \begin{cases} 0, & k < h,\\ a_{i-1}, & a_i \le k < a_{i+1},\ i \ge 2,\end{cases}
    \label{eq:hat_tau}
\end{equation}
and set
\begin{equation}
    \hat\tau_k \defeq k - d_k, \qquad w_k \defeq 1 + \hat\tau_k, \qquad \rho_k = c\,L\,w_k ,
    \label{eq:schedule_exact}
\end{equation}
so that $\hat\tau_k = k$ and $w_k = k+1$ for $k < h$, and $\hat\tau_k = \tau_k$ for $k \ge a_1$. The schedule is causal.

The \emph{delay window} of iteration $k'$ is $\{d_{k'}, \dots, k'-1\}$, the steps taken since its matrix was dispatched. The central technical tool is the Overlap Count, which counts how many later iterations have a given step $k$ in their delay window.

\begin{definition}[Overlap Count]
\label{def:overlap}
For $k \in \{0, \ldots, T-1\}$, set
$\Ov_k \defeq \sum_{k'=k+1}^{T-1} \mathbf{1}_{\{k' - \hat\tau_{k'} \le k\}}$.
\end{definition}

\begin{lemma}[Overlap Lemma]
\label{lem:overlap}
Assume the two-timeline model of \cref{def:timeline} and (A3). For every $k \ge a_1$, letting $i$ denote the curvature-job index with $a_i \le k < b_i$,
\begin{equation}
    \Ov_k \;\le\; \Delta_i + \Delta_{i+1} - 1 \;\le\; (1+\gamma)\Delta_i \;\le\; \gamma(1+\gamma)\,\tau_k.
\end{equation}
In particular $\Ov_k \le \Cov\, \tau_k$ with $\Cov \defeq \gamma(1+\gamma)$, and $\Cov \to 2$ as $\gamma \to 1$.
\end{lemma}

\begin{proof}
Fix $k$ with $a_i \le k < b_i$. For any $k' > k$ in curvature-job $j$ (so $a_j \le k' < b_j$), the curvature in use at $k'$ comes from job $j-1$, hence $k' - \tau_{k'} = a_{j-1}$. Thus the event $\{k' - \tau_{k'} \le k\} = \{a_{j-1} \le k\}$, equivalently $j - 1 \le i$, i.e.\ $j \le i+1$.

We count the $k' > k$ in jobs $j \le i+1$:
\begin{itemize}[itemsep=0.1em]
    \item \emph{Job $j = i$:} $k' \in (k, b_i) \cap \mathbb{N}$, so the count is at most $b_i - k - 1 \le b_i - a_i - 1 = \Delta_i - 1$.
    \item \emph{Job $j = i+1$:} $k' \in [b_i, b_{i+1}) \cap \mathbb{N}$, count is $\Delta_{i+1}$.
    \item \emph{Jobs $j < i$:} $k' < b_{i-1} \le a_i \le k$, contradicting $k' > k$. Zero count.
\end{itemize}
Summing, $\Ov_k \le \Delta_i + \Delta_{i+1} - 1$.

For the $\tau_k$ bound: $\tau_k = k - a_{i-1} \ge a_i - a_{i-1} = \Delta_{i-1}$. Applying (A3) at $i-1$ gives $\Delta_i \le \gamma\Delta_{i-1} \le \gamma\tau_k$, and applying it at $i$ gives $\Delta_{i+1} \le \gamma\Delta_i$:
\[
    \Ov_k \le \Delta_i + \Delta_{i+1} - 1 \le (1+\gamma)\Delta_i \le \gamma(1+\gamma)\tau_k. \qedhere
\]
\end{proof}

\paragraph{Startup and the weighted Hardy inequality.}
\cref{lem:overlap} is stated for $k \ge a_1$. Every iteration $k < h$ uses a matrix based at $x_0$ ($H_0$ for $k < a_1$, the job-$0$ output for $a_1 \le k < h$), so an early step lies in the delay window of every such iteration and \cref{lem:overlap} does not bound its charges. Instead, the prefix-sum operator is bounded in the weighted cubic norm of the analysis. With $r_k \defeq \norm{s_k}$ and the dispatch index $d_k$ of \eqref{eq:hat_tau}, write
\begin{equation}
    P_k \defeq \sum_{j=d_k}^{k-1} r_j, \qquad S_T \defeq \sum_{k<T} w_k r_k^3 ,
    \label{eq:Pk}
\end{equation}
with empty sums equal to zero. For $k < h$, $P_k$ is the full prefix sum.

\begin{lemma}[Weighted Hardy inequalities]
\label{lem:hardy}
Let $r_0, \dots, r_{n-1} \ge 0$ and $R_k = \sum_{j<k} r_j$. Then
\begin{equation}
    \begin{gathered}
    \sum_{k<n}\frac{R_k^3}{(k+1)^2} \;\le\; 27\sum_{k<n}(k+1)\,r_k^3,\\
    \sum_{k<n} R_k r_k^2 \;\le\; \sum_{k<n}(3k+1)\,r_k^3 .
    \end{gathered}
    \label{eq:hardy}
\end{equation}
No bound on $n$, on $R_k$ or on the individual steps is required.
\end{lemma}

\begin{lemma}[Tail overlap]
\label{lem:overlap_global}
Under (A3),
\begin{equation}
    \sum_{k=h}^{T-1}\ \sum_{j=d_k}^{k-1} r_j^3 \;\le\; \Cov\sum_{j=a_1}^{T-1} w_j\, r_j^3 ,
    \label{eq:tail_overlap}
\end{equation}
so no startup index $j < a_1$ is charged by any delay window with $k \ge h$.
\end{lemma}

\begin{proof}
For $k \ge h$ the window $\{d_k,\dots,k-1\}$ starts at $a_{i-1} \ge a_1$. A fixed $j \ge a_1$ in job $i \ge 1$ is charged at most $\Delta_i + \Delta_{i+1} - 1 \le \Cov\tau_j = \Cov\hat\tau_j$ times by the count of \cref{lem:overlap}, whose windows contain those used here.
\end{proof}

The proof of \cref{lem:hardy} is in \cref{app:init_phase}. While a matrix based at $x_0$ is in use, $\norm{\hess(x_k) - H_k} \le \bar\delta_0 + LP_k$, where $\bar\delta_0$ bounds the errors of $H_0$ and of the job-$0$ output. By \eqref{eq:hardy}, $\sum_{k<h}\big(\tfrac{L}{6}r_k^3 + \tfrac{L}{2}P_kr_k^2\big) \le \tfrac{3L}{2}\sum_{k<h}(k+1)r_k^3$, so $\rho_k = cL(k+1)$ pays for the whole prefix error without knowledge of the first job's duration, and $\bar\delta_0$ enters only through the effective gap
\begin{equation}
    \Ftilde_0 \;\defeq\; F_0 + \frac{25\pi^2\,\bar\delta_0^{\,3}}{81\,(cL)^2},
    \label{eq:Ftilde}
\end{equation}
which depends on $F_0$, $L$, $c$ and the initial curvature error only, and is uniform in the duration of the startup phase.

\subsection{Main convergence theorem}

\begin{theorem}[Exact asynchronous convergence]
\label{thm:exact}
Assume (A1)--(A5), exact gradients, and every published curvature exact ($E_i = 0$ for $i \ge 0$), with $H_0$ arbitrary subject to (A5); let $\Ftilde_0$ be given by \eqref{eq:Ftilde} with $\bar\delta_0 = \delta_0$. Run \cref{alg:split-client-cubic} with
\begin{equation}
    \rho_k = c\,L\,w_k, \qquad c \;\ge\; c_H(\gamma) \defeq \tfrac53(\Cov + 9),
    \label{eq:exact_schedule}
\end{equation}
with $w_k$ from \eqref{eq:schedule_exact}, and let $A(c,\Cov) \defeq \tfrac{3c}{20} - \tfrac32 - \tfrac{\Cov}{6}$, which is at least $c/20$ under \eqref{eq:exact_schedule}, and $K_H \defeq 27 + \Cov$.
\begin{enumerate}[itemsep=0.1em]
    \item \textbf{(First-order)} Sample $x_{\mathrm{out}} = x_{k+1}$ with probability $p_k \propto w_k^{-1/2}$, $k < T$, and let $\taubar = \tfrac1T\sum_{k<T}\hat\tau_k$. Then
    \[
        \E\!\left[\norm{\grad(x_{\mathrm{out}})}^{3/2}\right] \;\le\; \frac{C_0\, L^{1/2}\, \Ftilde_0\, \sqrt{1+\taubar}}{T},
    \]
    with $A = A(c,\Cov)$ and
    \[
        C_0 = \sqrt3\Big[\frac{(\tfrac{c+1}{2})^{3/2} + \sqrt{K_H}}{A} + c\sqrt{\frac{27}{50\,A}}\,\Big].
    \]
    \item \textbf{(Second-order, after startup)} If $T > h$, let $N = T - h$ and $W_+ = 1 + \max_{h \le k < T}\tau_k$. There exists $k^\star \in \{h, \dots, T-1\}$ such that
    \[
        \lmin(\hess(x_{k^\star + 1})) \;\ge\; -\,C_2\,\frac{L^{2/3}\, \Ftilde_0^{1/3}\, W_+^{2/3}}{N^{1/3}},
    \]
    with $C_2 = A(c,\Cov)^{-1/3}\big(\tfrac{c}{2}+\Cov^{1/3}\big)$.
\end{enumerate}
Both constants depend on $c$ and on $\Cov = \gamma(1+\gamma)$; neither can be taken absolute (\cref{rem:C2}). At $c = c_H$ one has $A = c/20 = (\Cov+9)/12$, $C_0 = \mathcal{O}(\sqrt{\Cov})$ and $C_2 = \mathcal{O}(\Cov^{2/3})$; for $\gamma = 1$ ($\Cov = 2$) they equal $91.3$ and $10.7$, and for $\gamma = 1.1$ ($\Cov = 2.31$) they equal $92.1$ and $11.0$. The delay exponent in the second-order bound is $2/3$.
\end{theorem}

The gap $\Ftilde_0$ is uniform in the startup duration, and the algorithm needs neither the first job's duration nor a positive semidefinite $H_0$; the price is the larger multiplier $c_H$. Restricting the second-order index to $k \ge h$ changes no update; an all-horizon variant with an additive startup term is in \cref{app:proof_exact_hess}.

For fixed $\gamma$, the first-order horizon factor depends on the realized average age rather than the maximum age. The constants retain their dependence on $\gamma$, so a rare large increase in job duration enlarges them. \cref{app:avg_vs_max} illustrates the contrast with a $\tmax$-based analysis.

\subsection{Wall-clock comparison}
\label{sec:wall_clock_comparison}

\paragraph{Job cost and curvature age are different quantities.} Recall that $D$ is the duration of a curvature job in iteration units and $\tau_k$ the \emph{age} of the curvature in use. If every job takes $D$ units, the steady-state ages within a job run $D, D+1, \dots, 2D-1$, so the time-averaged age is $(3D-1)/2$ while the job cost is $D$. For variable durations, a complete job-$i$ block contributes $\sum_{k=a_i}^{a_{i+1}-1}\tau_k = \Delta_i\Delta_{i-1} + \tfrac{\Delta_i(\Delta_i-1)}{2}$, so the average age is length-weighted and involves products and second moments of durations.

\paragraph{Fixed-cost comparison.} Under a stipulated common iteration cost and a fixed job duration $D$, the startup ages are at most $D-1$ and all later ages at most $2D-1$, so $1 + \taubar \le 2D$. \Cref{thm:exact} then gives, with $B_A \defeq C_0 L^{1/2}\Ftilde_0$,
\begin{equation}
    T_A \ge B_A\sqrt{2D}\,\eps^{-3/2} \;\Longrightarrow\; \E\norm{\grad(x_{\mathrm{out}})}^{3/2} \le \eps^{3/2},
    \label{eq:speedup}
\end{equation}
where $\Ftilde_0$ does not depend on $D$, so the startup phase adds no hidden factor, provided $H_0$ is cheap to apply (e.g.\ a multiple of the identity); a dense $H_0$ adds its own setup cost. With $B_V$ and $B_L$ the corresponding iteration-complexity constants of Vanilla and Lazy cubic Newton, the leading $\eps^{-3/2}$ terms for fixed job duration $D$ are:
\begin{center}\small
\begin{tabular}{lc}
\toprule
Method & Leading term \\
\midrule
Vanilla Cubic & $B_V(1+D)\,\eps^{-3/2}$ \\
Lazy Hessian (reuse $m$) & $B_L(\sqrt{m} + D/\sqrt{m})\,\eps^{-3/2}$ \\
Lazy (optimized $m$) & $\mathcal{O}(B_L\sqrt{D}\,\eps^{-3/2})$ \\
\textbf{Split-Client} & $\mathcal{O}(B_A\sqrt{D}\,\eps^{-3/2})$ \\
\bottomrule
\end{tabular}
\end{center}
The Lazy iteration complexity is $\mathcal{O}(\sqrt{m}\,\eps^{-3/2})$ \citep{doikov2023lazy}, and minimizing $\sqrt{m} + D/\sqrt{m}$ gives $m$ near $D$.

\begin{corollary}[Comparison with the arithmetic reuse period]
\label{cor:lazy_md}
In the fixed-cost model with job duration $D$, let Lazy Hessian use $m = d$, the reuse period that is optimal when a curvature job costs $d$ gradient evaluations \citep{doikov2023lazy}. The ratio of its leading term to that of Split-Client in \eqref{eq:speedup} is
\[
\frac{B_L\big(\sqrt{d} + D/\sqrt{d}\big)}{B_A\sqrt{2D}} \;=\; \frac{B_L}{\sqrt2\,B_A}\Big(\sqrt{\tfrac{d}{D}} + \sqrt{\tfrac{D}{d}}\Big).
\]
\end{corollary}

The bracket equals $2$ when $D = d$ and grows like $\sqrt{\max(D/d,\, d/D)}$ otherwise. In the dense arithmetic model a curvature job costs $\Theta(nd^2 + d^3)$ and an iteration $\Theta(nd + d^2)$, so $D = \Theta(d)$ and the factor is a constant; it reflects how far the measured $D$ is from $d$, to which Split-Client adapts without knowing $D$. With $m$ tuned to $D$, the two leading terms have the same order. These are sufficient upper bounds with method-specific constants, not realized runtimes, and they suppress initialization and integer rounding (a finite Lazy run of $K$ iterations pays $D\lceil K/m\rceil$ for curvature). We claim no advantage over variable-reuse policies under variable job costs. \cref{app:phys_time} gives an elapsed-time bound for arbitrary iteration costs.

\section{Theory: Inexact Curvature}
\label{sec:theory_inexact}
% ============================================================================

We now allow inexact curvature, $\Ht_k = \hess(x_{d_k}) + E_i$ for $k$ in job $i+1$ ($i \ge 0$), and an inexact gradient $\gt_k$. Two error models arise:
\begin{itemize}[itemsep=0.1em]
    \item \textbf{Structured (causal adaptive):} the curvature error of job $i$ is at most $\alpha\norm{s_{a_i - 1}}$, a step known when the job is dispatched. This covers adaptively subsampled Hessians (\S\ref{sec:subsampled}). Memory-based Quasi-Newton matrices have a related, directional error structure, recorded in \cref{app:qn_full}.
    \item \textbf{Uniform (worst-case):} $\norm{E_i} \le \delta$ for a constant $\delta$. Robustness fallback (\S\ref{sec:robust}).
\end{itemize}
For bounded average age and fixed problem constants, the stationarity measure $\E\norm{\grad(x_{\mathrm{out}})}^{3/2}$ decays as $\mathcal{O}(T^{-1})$ under the structured model (\cref{thm:unified}) and $\widetilde{\mathcal{O}}(T^{-3/4})$ under the uniform one (\cref{thm:robust}).

\begin{figure*}[t]
  \centering
  \begin{minipage}[b]{.32\linewidth}
    \centering
    \includegraphics[width=\linewidth]{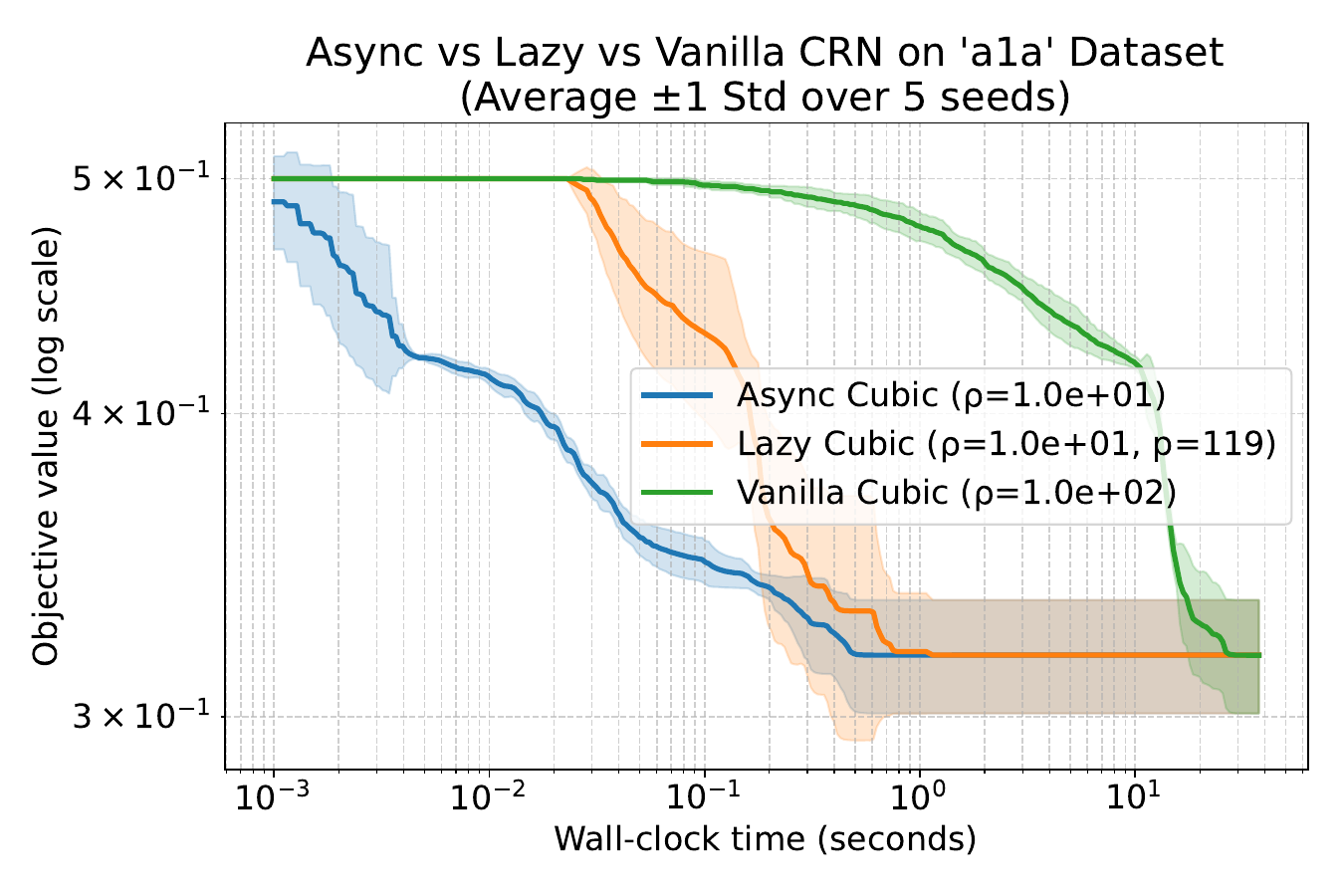}
    \caption{Geman--McClure (\texttt{a1a}, $d{=}119$): objective versus wall-clock time.}
    \label{fig:nonconvex-loss}
  \end{minipage}\hfill
  \begin{minipage}[b]{.32\linewidth}
    \centering
    \includegraphics[width=\linewidth]{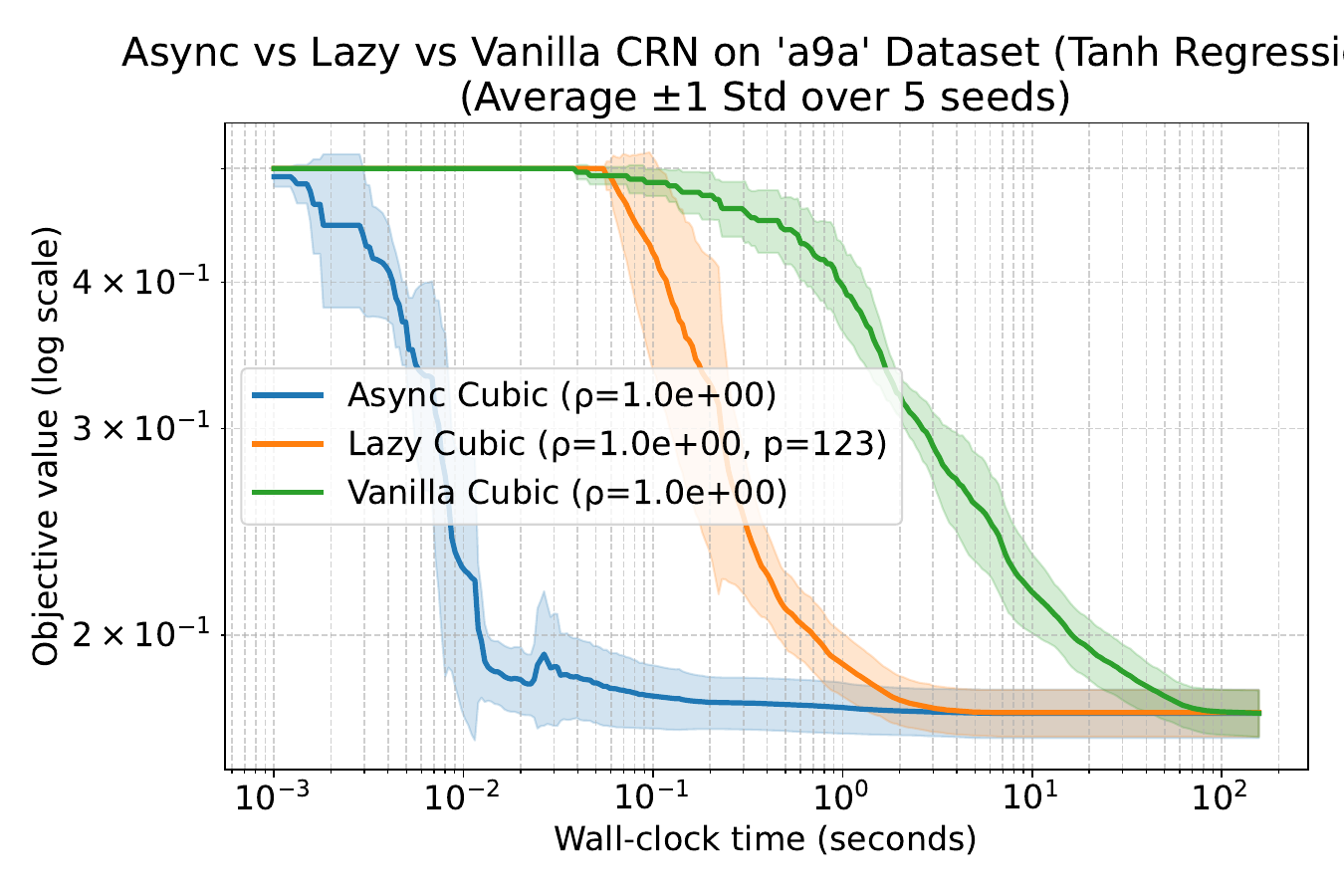}
    \caption{Tanh (\texttt{a9a}, $d{=}123$): objective versus wall-clock time.}
    \label{fig:tanh-loss}
  \end{minipage}\hfill
  \begin{minipage}[b]{.32\linewidth}
    \centering
    \includegraphics[width=\linewidth]{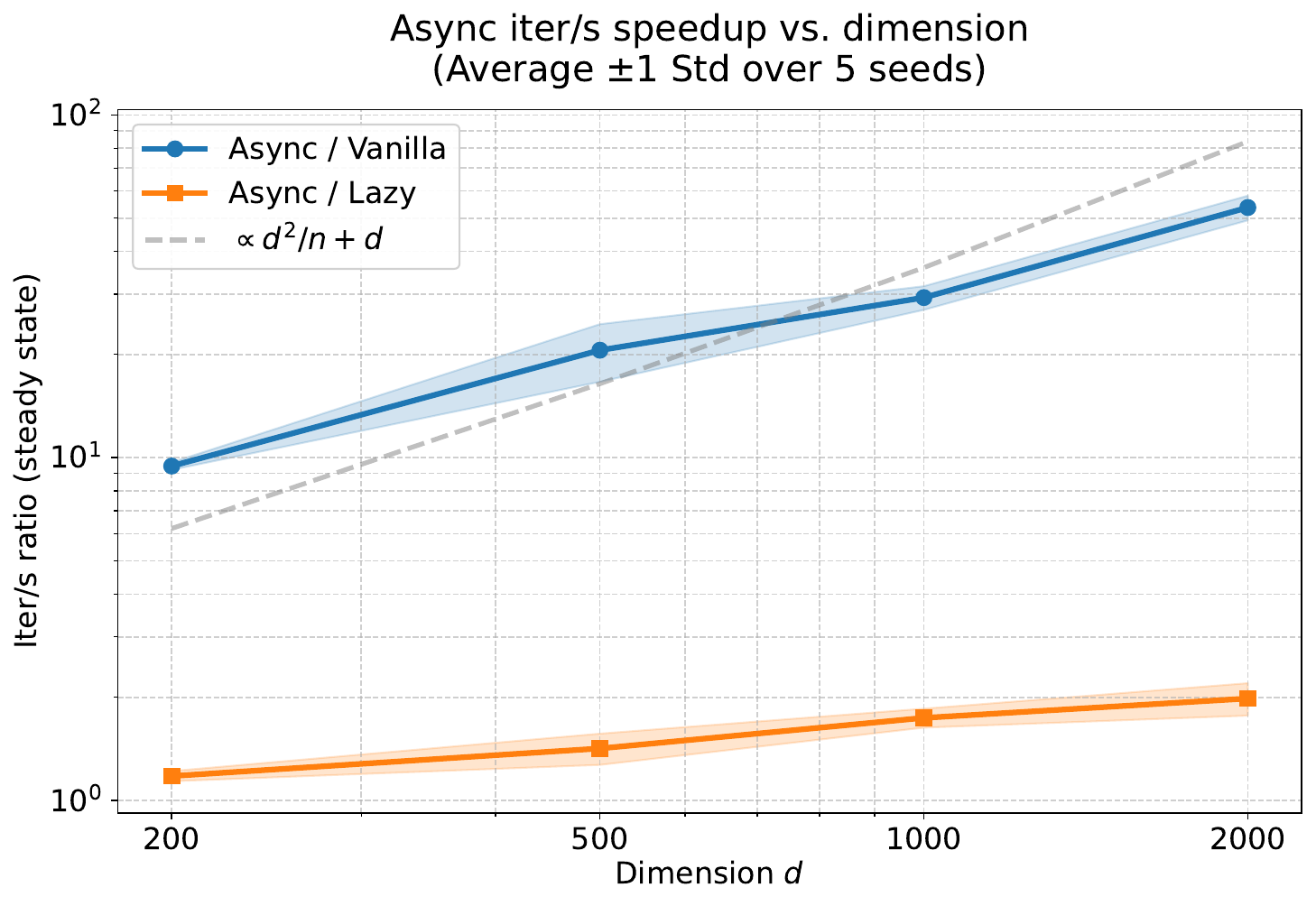}
    \caption{Iteration-throughput ratios of Split-Client over Vanilla and Lazy vs.\ $d$ ($n{=}5000$, 5 seeds); not time-to-target ratios.}
    \label{fig:speedup-vs-d}
  \end{minipage}
\end{figure*}

\subsection{Structured inexactness: unified theorem}
\label{sec:structured}

\begin{enumerate}[label=(A\arabic*),start=6,itemsep=0em]
    \item \textbf{Causal adaptive inexactness:} there exist $\alpha, \beta \ge 0$ such that
      \begin{align*}
         &\norm{E_i} \le \alpha\, r_{a_i - 1} \quad (i \ge 1), \\
         &\norm{\gt_k - \grad(x_k)} \le \beta\, r_{k-1}^2 \quad (k \ge 1), \\
         &\gt_0 = \grad(x_0).
      \end{align*}
\end{enumerate}

Both tolerances are known before sampling (the Hessian one uses the step preceding dispatch, the gradient one the previous step), as in the previous-step conditions of \citet{wang2018inexact} and unlike conditions on the step being computed \citep{kohler2017sub}; \citet{xu2020newton} enforce a fixed accuracy a priori. Our analysis also accounts for reusing a matrix computed at a dispatch snapshot. The convention $r_{-1} = 0$ records that the first gradient is exact.

\paragraph{Realizability of (A6).}
For finite-sum objectives, (A6) is realized by adaptive subsampling with these tolerances: under bounded-deviation assumptions on the components, matrix concentration gives batch sizes for which each bound holds with high probability, conditionally on the past, and a union bound gives a simultaneous good event (\cref{sec:subsampled,app:a6_realize}).

\begin{theorem}[Unified structured convergence]
\label{thm:unified}
Assume (A1)--(A5), that the two matrices based at $x_0$ (the initial matrix $H_0$ and the output of job $0$) have error at most $\bar\delta_0$, and that (A6) holds throughout the run. Run \cref{alg:split-client-cubic} with the sampled gradient and
\begin{equation}
    \rho_k = c\,L\,w_k + (4+\gamma^2)\alpha + 2\gamma^2\alpha\hat\tau_k + 6\beta, \qquad c \ge c_H(\gamma),
    \label{eq:unified_schedule}
\end{equation}
and sample $x_{\mathrm{out}} = x_{k+1}$ with probability $p_k \propto \rho_k^{-1/2}$. Let $\bar\rho \defeq \tfrac1T\sum_{k<T}\rho_k = cL(1+\taubar) + (4+\gamma^2)\alpha + 2\gamma^2\alpha\taubar + 6\beta$, $K_H = 27 + \Cov$, and $\Ftilde_0$ as in \eqref{eq:Ftilde}. Then
\[
    \E\!\left[\norm{\grad(x_{\mathrm{out}})}^{3/2}\right] \;\le\; \frac{C_1\, \Ftilde_0\, \sqrt{\bar\rho}}{T},
\]
\[
    C_1 = 2\,\Big[20\Big(\tfrac{17}{16} + \tfrac{\gamma^2}{8} + \tfrac{\sqrt{K_H}}{c^{3/2}} + \tfrac{1}{6\sqrt6}\Big) + \sqrt{\tfrac{54}{5}}\,\Big].
\]
If moreover $T > h$, with $N = T - h$, $W_+ = 1 + \max_{h\le k<T}\tau_k$ and $\rho_{+} = \max_{h \le k < T}\rho_k$, some $k^\star \in \{h,\dots,T-1\}$ satisfies
\[
    -\lmin(\hess(x_{k^\star + 1})) \;\le\; \frac{(20\Ftilde_0)^{1/3}}{N^{1/3}}\,\Phi,
\]
with $\Phi = \tfrac12\rho_+^{2/3} + \alpha\big(\tfrac{\gamma^2}{cL}\big)^{1/3} + L^{2/3}W_+^{2/3}\big(\tfrac{\Cov}{c}\big)^{1/3}$.
Both coefficients depend on $c$, $\gamma$ and $\Cov$. The first-order horizon factor involves the average age: the linear-in-$\hat\tau_k$ inflation of \eqref{eq:unified_schedule}, needed because the error of each matrix is charged to its reference step up to $\gamma^2$ times the weight of that step (\cref{lem:reuse}), is averaged to $\taubar$ by the sampling $p_k \propto \rho_k^{-1/2}$. With $\alpha = \beta = 0$ the sharper constants of \cref{thm:exact} apply.

These are deterministic statements, conditional on the run, that hold whenever (A6) holds at every iteration. When (A6) is realized by subsampling (\cref{sec:subsampled}), it holds simultaneously with probability at least $1 - \eta$. \cref{app:a6_realize} gives an unconditional expectation bound under an additional gradient bound covering failed-oracle trajectories; without that condition, the guarantee holds on the simultaneous good event. Proof in \cref{app:proof_unified}.
\end{theorem}

\subsection{Extension: Subsampled Cubic Newton}
\label{sec:subsampled}
Let $f = \tfrac1n\sum_{j=1}^nf_j$ with $\norm{\nabla f_j(x) - \grad(x)} \le \sigma$ and $\norm{\nabla^2 f_j(x) - \hess(x)} \le \varsigma$ at every sampled point. For a tolerance $t > 0$ and failure probability $p$, draw fresh independent uniform samples with replacement, of count
\begin{equation}
    \begin{aligned}
    b_H(t,p) &= \min\Big\{n,\ \max\Big\{1, \Big\lceil\tfrac{8\varsigma^2}{t^2}\log\tfrac{2d}{p}\Big\rceil\Big\}\Big\},\\
    b_g(t,p) &= \min\Big\{n,\ \max\Big\{1, \Big\lceil\tfrac{8\sigma^2}{t^2}\log\tfrac{2(d+1)}{p}\Big\rceil\Big\}\Big\},
    \end{aligned}
    \label{eq:batch_sizes}
\end{equation}
and evaluate the full finite sum exactly when the count reaches $n$ or when $t = 0$. The tolerances are $t_H = \alpha r_{a_i - 1}$ for job $i \ge 1$ and $t_g = \beta r_{k-1}^2$ for iteration $k \ge 1$; the first gradient and job $0$ are computed exactly, so $\bar\delta_0 = \delta_0$, and job $0$ runs in the background without delaying the first steps. Matrix Hoeffding \citep[Theorem~1.3]{tropp2012user} bounds the spectral error of a Hessian batch, and the same inequality applied to the symmetric dilation of a gradient sample, in dimension $d+1$, bounds the gradient error; this gives \eqref{eq:batch_sizes}. Every tolerance and count is determined before its batch is drawn, so the concentration bounds hold conditionally at adaptive iterates, and (A6) holds at each call with probability at least $1-p$. Sampling does not by itself establish (A3), which remains a separate requirement on the resulting schedule.

\begin{proposition}[Aggregate sample budget]
\label{prop:budget}
Count the curvature jobs dispatched at indices $a_i < T$, job $0$ included; positive integer durations give at most $N_H \le T$ of them. Each gradient or curvature call uses at most $n$ component evaluations, so
\begin{equation}
    \mathcal{B}^g_T \;\le\; n\,T, \qquad \mathcal{B}^H_T \;\le\; n\,N_H \;\le\; n\,T .
    \label{eq:budget}
\end{equation}
\end{proposition}

\subsection{Robustness fallback: uniform bounded error}
\label{sec:robust}

When (A6) is unavailable, an inflated regularization schedule gives the following guarantee under a uniform curvature error:

\begin{theorem}[Anytime robustness]
\label{thm:robust}
Assume (A1)--(A4), exact gradients, and a uniform curvature error: $\norm{H_0 - \hess(x_0)} \le \delta$ and $\norm{E_i} \le \delta$ for every $i \ge 0$. With $\rho_k = \max(c\,L\,w_k, \nu\delta^{3/2}\sqrt{k+1})$ for $c \ge c_H(\gamma)$ and $\nu > 0$, and $p_k \propto w_k^{-1/2}$, define the noise-inclusive work budget
\begin{equation}
    \mathcal{K}(T) \;\defeq\; 20\,\Ftilde_0 \;+\; \frac{1000\,(1+\ln T)}{27\,\nu^2}.
    \label{eq:KT}
\end{equation}
Then
\[
    \begin{aligned}
    \frac{\E\!\left[\norm{\grad(x_{\mathrm{out}})}^{3/2}\right]}{\sqrt3\,\sqrt{1+\taubar}} \;\le\;& \Big(\tfrac65\sqrt c + \tfrac{\sqrt{K_H}}{c}\Big)\frac{L^{1/2}\,\mathcal{K}(T)}{T}\\
    &+ \delta^{3/4}\,\frac{\sqrt\nu\,\mathcal{K}(T) + \sqrt{2\mathcal{K}(T)/\nu}}{T^{3/4}},
    \end{aligned}
\]
with $\Ftilde_0$ from \eqref{eq:Ftilde} for $\bar\delta_0 = \delta$. For bounded average age this is $\widetilde{\mathcal{O}}(T^{-3/4})$; a fixed-horizon variant with $\sqrt{T}$ in place of $\sqrt{k+1}$ removes the logarithm. The growing noise budget is needed in this analysis: for $f \equiv 0$, unit job durations and every matrix equal to $-\delta I$, each step has norm $2\delta/\rho_k$, so $\sum_k \rho_kr_k^3$ grows like $8\ln T/\nu^2$ while $\Ftilde_0$ stays fixed. The schedule requires a bound $\delta$ on the noise. Proof in \cref{app:proof_robust}.
\end{theorem}

% ============================================================================
\section{Experiments}
\label{sec:experiments}
% ============================================================================

We evaluate Split-Client against Vanilla and Lazy Cubic Newton \citep{doikov2023lazy} on two regression problems with nonconvex components (a Geman--McClure penalty and a tanh model) where factorization is the wall-clock bottleneck, using Python \texttt{multiprocessing} with shared-memory handoff. Each method tunes $\rho$ on a log-grid, and Lazy reuses each Hessian for $m = d$ steps. Split-Client is labeled ``Async'' in the figures. Full setup in \cref{app:additional_experiments}.

\paragraph{Wall-clock convergence.} On the \texttt{a1a} ($d{=}119$) and \texttt{a9a} ($d{=}123$) LIBSVM datasets \citep{chang2011libsvm} we run Geman--McClure regression (\cref{fig:nonconvex-loss}) and tanh regression (\cref{fig:tanh-loss}), using randomized SVD for indefinite Hessians. On both datasets Async decreases the objective earliest. Its lead over Lazy is largest at intermediate objective values and shrinks as all three methods approach a common final plateau, which Vanilla reaches last. Objective curves do not certify a stationary point, and the runs use tuned constant regularization and low-rank curvature, so they evaluate an implementation of the architecture rather than the exact schedules of the theorems.

\paragraph{Diagnostics.} On these runs, ages observed at curvature-update times lie in $\{2,\dots,9\}$, and peak memory is about twice Vanilla's, as double buffering predicts (\cref{app:assumption_diagnostics}).

\paragraph{Speedup grows with dimension.}
Against a gradient pass a curvature job costs $q_{\mathrm{grad}} = \Theta(d + d^2/n)$ (\cref{app:phys_time}), the denominator of the reference curve in \cref{fig:speedup-vs-d}; the Hessian is assembled in full at every job, without caching $A^\top A/n$. From $d{=}200$ to $2000$ the throughput ratio grows from $9.5\times$ to $53.7\times$ over Vanilla and from about $1.2\times$ to $2\times$ over Lazy; this tests the cost model, not the rate.

% ============================================================================
\section{Discussion and Conclusion}
\label{sec:discussion}

Split-Client separates curvature computation from gradient steps, so that the cost of forming and factorizing the Hessian appears as the age of the curvature in use rather than as blocking time. A regularization proportional to the observed age, an overlap count and a weighted Hardy inequality give a first-order bound in terms of the average age for any bounded-error symmetric initial matrix and, under a common fixed-cost model, the curvature-cost exponent of optimally tuned reuse without a reuse period chosen in advance.

\paragraph{Limitations.} The schedules depend on the Hessian Lipschitz constant $L$, which is rarely known; the experiments therefore tune the regularization, and fully adaptive schedules would remove this step. A Lipschitz Hessian is standard for cubic regularization; relaxing it, for instance to H\"older-continuous Hessians \citep{grapiglia2021regularized}, would broaden the results. The second-order bound holds at an index selected by the analysis, not at the sampled output, and the comparison with Lazy Hessian concerns leading terms of upper bounds.

\paragraph{Open directions.} The method factorizes a dense Hessian, so all claims concern moderate $d$. Matrix-free curvature built from Hessian-vector products, or structured preconditioners such as K-FAC \citep{martens2015optimizing} and Shampoo \citep{gupta2018shampoo}, would need a new cost model and possibly a different notion of age. Other directions are several curvature jobs in flight, rates exploiting the directional structure of Quasi-Newton matrices (\cref{app:qn_full}), and experiments at larger $d$.

\subsection*{AI use statement}
We used a large language model as an assistant for restructuring and editing the manuscript, checking and tightening proofs, drafting parts of the experiment code, and suggesting related work. The authors verified all such output. The research ideas, the method and the main results are the authors' own, and the authors take full responsibility for the content of this work.

\bibliographystyle{plainnat}
\bibliography{iclr2026_conference}

\section*{Checklist}

\begin{enumerate}

  \item For all models and algorithms presented, check if you include:
  \begin{enumerate}
    \item A clear description of the mathematical setting, assumptions, algorithm, and/or model. [Yes] The two-timeline model is \cref{def:timeline}, the assumptions are (A1)--(A6) in \cref{sec:theory_exact,sec:theory_inexact}, and the pseudocode is \cref{alg:split-client-cubic}.
    \item An analysis of the properties and complexity (time, space, sample size) of any algorithm. [Yes] Convergence rates are \cref{thm:exact,thm:unified,thm:robust}; wall-clock scaling is in \cref{sec:wall_clock_comparison}; the sample budget is \cref{prop:budget}; memory is discussed in \cref{app:assumption_diagnostics}.
    \item (Optional) Anonymized source code, with specification of all dependencies, including external libraries. [No] The code will be released after the final decision.
  \end{enumerate}

  \item For any theoretical claim, check if you include:
  \begin{enumerate}
    \item Statements of the full set of assumptions of all theoretical results. [Yes] Each theorem lists its assumptions.
    \item Complete proofs of all theoretical results. [Yes] \cref{app:proof_exact,app:proof_unified,app:budget_proof,app:proof_robust}.
    \item Clear explanations of any assumptions. [Yes] \cref{sec:theory_exact,sec:theory_inexact}.
  \end{enumerate}

  \item For all figures and tables that present empirical results, check if you include:
  \begin{enumerate}
    \item The code, data, and instructions needed to reproduce the main experimental results (either in the supplemental material or as a URL). [No] The datasets are public LIBSVM datasets, and the code will be released after the final decision.
    \item All the training details (e.g., data splits, hyperparameters, how they were chosen). [Yes] \Cref{app:additional_experiments} lists the selected hyperparameters and the tuning protocol: $\rho$ was searched over powers of $10$ and then refined around the best value.
    \item A clear definition of the specific measure or statistics and error bars (e.g., with respect to the random seed after running experiments multiple times). [Yes] Curves average five seeds with $\pm1$ standard deviation bands.
    \item A description of the computing infrastructure used. (e.g., type of GPUs, internal cluster, or cloud provider). [No] All experiments were run on Google Colab (\cref{app:additional_experiments}); the exact CPU model and thread settings are not reported.
  \end{enumerate}

  \item If you are using existing assets (e.g., code, data, models) or curating/releasing new assets, check if you include:
  \begin{enumerate}
    \item Citations of the creator If your work uses existing assets. [Yes] LIBSVM datasets are cited.
    \item The license information of the assets, if applicable. [No] The LIBSVM datasets are publicly available and cited; their license is not stated in the paper.
    \item New assets either in the supplemental material or as a URL, if applicable. [Not Applicable]
    \item Information about consent from data providers/curators. [Not Applicable]
    \item Discussion of sensible content if applicable, e.g., personally identifiable information or offensive content. [Not Applicable]
  \end{enumerate}

  \item If you used crowdsourcing or conducted research with human subjects, check if you include:
  \begin{enumerate}
    \item The full text of instructions given to participants and screenshots. [Not Applicable]
    \item Descriptions of potential participant risks, with links to Institutional Review Board (IRB) approvals if applicable. [Not Applicable]
    \item The estimated hourly wage paid to participants and the total amount spent on participant compensation. [Not Applicable]
  \end{enumerate}

\end{enumerate}

% ============================================================================

%========================================

%========================================
\clearpage
\appendix
\thispagestyle{empty}
\onecolumn
\aistatstitle{Second-Order Optimization with Asynchronous Curvature:\\ Supplementary Materials}

% ============================================================================
\section{Algorithm pseudocode}
\label{app:algorithm}
% ============================================================================

\begin{algorithm}[H]
\caption{Split-Client Cubic Regularization}
\label{alg:split-client-cubic}
\begin{algorithmic}[1]
\STATE \textbf{Input:} initial point $x_0$, initial matrix $H_0$ (e.g.\ a multiple of the identity), multiplier $c$, Lipschitz constant $L$ (and $\alpha,\beta$ or $\nu,\delta$ for the inexact variants).
\STATE \textbf{Curvature client} (concurrent): \textbf{loop:} read the current iterate $x_{k}$ and its index $k$; compute and factorize the curvature at $x_{k}$ (in the subsampled variant with $b_H(\alpha r_{k-1}, p)$ samples, and exactly for the first job); publish it with dispatch index $k$.
\STATE \textbf{Gradient client:} \textbf{for} $k = 0, 1, 2, \ldots$ \textbf{do}
\STATE \quad Receive $g_k$: the exact gradient, or in the subsampled variant a batch of size $b_g(\beta r_{k-1}^2, p)$ from \eqref{eq:batch_sizes}, with $g_0$ exact.
\STATE \quad If a new curvature has been published, load it and its dispatch index $d$.
\STATE \quad Set the age $\hat\tau_k = k$ while $H_0$ is in use, and $\hat\tau_k = k - d$ afterwards.
\STATE \quad Set $\rho_k$ by the schedule of the relevant theorem, evaluated at $\hat\tau_k$.
\STATE \quad Compute $s_k \in \arg\min_s\big\{\inner{g_k}{s} + \tfrac12\inner{H_ks}{s} + \tfrac{\rho_k}{6}\norm{s}^3\big\}$, with $g_k$ the gradient received above and $H_k$ the (symmetric) matrix in use; set $x_{k+1} \leftarrow x_k + s_k$.
\STATE \textbf{Output:} $x_{k+1}$ with $k$ drawn from the sampling distribution of the relevant theorem.
\end{algorithmic}
\end{algorithm}

% ============================================================================
\section{Proof of Theorem \ref{thm:exact} (Exact Delayed Convergence)}
\label{app:proof_exact}
% ============================================================================

The proof has four parts: (\ref{app:proof_exact_descent}) a per-iteration descent inequality, (\ref{app:init_phase}) the weighted Hardy inequalities controlling startup, (\ref{app:proof_exact_work}) a global work bound, and (\ref{app:proof_exact_grad})--(\ref{app:proof_exact_hess}) the first- and second-order rates.

\paragraph{Conventions.} Throughout, $H_k$ is the matrix used at iteration $k$, $h = a_2$, and $d_k$, $\hat\tau_k$, $w_k$, $P_k$, $S_T$ are as in \eqref{eq:hat_tau}--\eqref{eq:Pk}; in every schedule, weight and average, $\tau_k$ stands for $\hat\tau_k$. Every step is a global minimizer of the cubic model, so
\begin{equation}
    \gt_k + H_k s_k + \tfrac{\rho_k}{2}r_k s_k = 0, \qquad H_k + \tfrac{\rho_k}{2}r_k I \succeq 0 ,
    \label{eq:cubic_certificates}
\end{equation}
with $\gt_k = \grad(x_k)$ in this appendix. The second condition is a property of global cubic minimizers and is not implied by second-order local optimality. For $k < h$ the matrix is based at $x_0$ and $\norm{\hess(x_k) - H_k} \le \bar\delta_0 + L\norm{x_k - x_0} \le \bar\delta_0 + LP_k$; for $k \ge h$ it is the exact Hessian at $x_{d_k}$ and $\norm{\hess(x_k) - H_k} \le LP_k$. Hence, for every $k$,
\begin{equation}
    \norm{\hess(x_k) - H_k} \;\le\; LP_k + \bar\delta_0\,\mathbf{1}_{k<h}.
    \label{eq:error_bound}
\end{equation}

\subsection{Per-iteration descent}
\label{app:proof_exact_descent}

\begin{lemma}[Descent inequality]
\label{lem:descent}
Under (A1), for every $k$,
\begin{equation}
f(x_{k+1}) - f(x_k) \;\le\; -\Big(\tfrac{\rho_k}{4} - \tfrac{L}{6}\Big)r_k^3 + \tfrac12\norm{\hess(x_k) - H_k}\,r_k^2 .
\label{eq:descent-main}
\end{equation}
\end{lemma}

\begin{proof}
Smoothness gives $f(x_{k+1}) \le f(x_k) + \inner{\grad(x_k)}{s_k} + \tfrac12\inner{\hess(x_k)s_k}{s_k} + \tfrac{L}{6}r_k^3$. The first condition in \eqref{eq:cubic_certificates} gives $\inner{\grad(x_k)}{s_k} = -\inner{H_ks_k}{s_k} - \tfrac{\rho_k}{2}r_k^3$, so
\begin{equation}
f(x_{k+1}) - f(x_k) \;\le\; \tfrac{1}{2}\inner{(\hess(x_k) - H_k) s_k}{s_k} - \tfrac12\inner{H_k s_k}{s_k} + \tfrac{L}{6}r_k^3 - \tfrac{\rho_k}{2}r_k^3 .
\label{eq:desc-raw}
\end{equation}
The term $-\tfrac12\inner{H_ks_k}{s_k}$ may be positive, since $H_k$ may be indefinite. The second condition in \eqref{eq:cubic_certificates} gives
\begin{equation}
\inner{H_k s_k}{s_k} \;\ge\; -\tfrac{\rho_k}{2}r_k^3 ,
\label{eq:cubic_second_order}
\end{equation}
whence
\begin{equation}
f(x_{k+1}) - f(x_k) \;\le\; \tfrac{1}{2}\inner{(\hess(x_k) - H_k) s_k}{s_k} + \Big(\tfrac{L}{6} - \tfrac{\rho_k}{4}\Big)r_k^3 ,
\label{eq:desc-intermediate}
\end{equation}
which implies \eqref{eq:descent-main}.
\end{proof}

For $k \ge h$, Young's inequality $ab^2 \le \tfrac13a^3 + \tfrac23b^3$ applied to each summand of $\tfrac{L}{2}P_kr_k^2$ gives
\begin{equation}
\tfrac{L}{2}P_k r_k^2 \;\le\; \tfrac{L}{6}\sum_{j=d_k}^{k-1}r_j^3 + \tfrac{L\hat\tau_k}{3}r_k^3 ,
\label{eq:young_tail}
\end{equation}
and together with $\tfrac{L}{6}r_k^3$ the current-step coefficient is at most $\tfrac{L}{3}w_k$, since $\tfrac16 + \tfrac{\hat\tau_k}{3} \le \tfrac13(1+\hat\tau_k)$.

\subsection{Startup: the weighted Hardy inequalities}
\label{app:init_phase}

\begin{proof}[Proof of \cref{lem:hardy}]
Let $r(t) = r_k$ for $t \in [k, k+1)$ and $R(t) = \int_0^t r(u)\,du$, and put
\[
I = \int_0^n \frac{R(t)^3}{t^2}\,dt, \qquad J = \int_0^n t\,r(t)^3\,dt = \sum_{k<n}\big(k+\tfrac12\big)r_k^3 .
\]
The boundary term at zero vanishes since $R(t) = r_0t$ near $0$. Integration by parts and H\"older's inequality give
\[
I = -\frac{R(n)^3}{n} + 3\int_0^n\frac{R(t)^2r(t)}{t}\,dt \;\le\; 3\,I^{2/3}J^{1/3},
\]
so $I \le 27J$ (the case $I = 0$ is immediate). For $k \ge 1$, $R(t) \ge R_k$ on $[k, k+1]$ and $\int_k^{k+1}t^{-2}\,dt = \tfrac{1}{k(k+1)} \ge \tfrac{1}{(k+1)^2}$, while the $k = 0$ summand vanishes; this proves the first inequality in \eqref{eq:hardy}. H\"older's inequality also gives $\int_0^n R(t)r(t)^2\,dt \le I^{1/3}J^{2/3} \le 3J$, and $\int_0^n R(t)r(t)^2\,dt = \sum_k R_kr_k^2 + \tfrac12\sum_k r_k^3$. Hence $\sum_kR_kr_k^2 \le \sum_k(3k + \tfrac32 - \tfrac12)r_k^3$, the second inequality.
\end{proof}

\begin{lemma}[Startup phase]
\label{lem:init}
Let $n \le h$. Under (A1) and $\rho_k = cL(k+1)$ for $k < h$,
\begin{equation}
\sum_{k<n}\Big(\tfrac{L}{6}r_k^3 + \tfrac{L}{2}P_kr_k^2\Big) \;\le\; \tfrac{3L}{2}\sum_{k<n}w_kr_k^3,
\qquad
\sum_{k<n}\Big(\tfrac{\bar\delta_0}{2}r_k^2 - \tfrac{cL}{10}w_kr_k^3\Big) \;\le\; \frac{25\pi^2\,\bar\delta_0^{\,3}}{81\,(cL)^2}.
\label{eq:descent-startup}
\end{equation}
\end{lemma}

\begin{proof}
For $k < h$, $P_k = R_k$ and $w_k = k+1$. The second inequality in \eqref{eq:hardy} gives $\sum_{k<n}\big(\tfrac{L}{6}r_k^3 + \tfrac{L}{2}R_kr_k^2\big) \le L\sum_{k<n}\big(\tfrac32k + \tfrac23\big)r_k^3 \le \tfrac{3L}{2}\sum_{k<n}(k+1)r_k^3$. For the second, $\max_{u\ge0}(Bu^2 - Au^3) = \tfrac{4B^3}{27A^2}$ with $B = \bar\delta_0/2$ and $A = cL(k+1)/10$ gives $\tfrac{50\bar\delta_0^3}{27(cL)^2(k+1)^2}$, and $\sum_{k\ge0}(k+1)^{-2} = \pi^2/6$ gives the total $\tfrac{50\pi^2}{162}\tfrac{\bar\delta_0^3}{(cL)^2} = \tfrac{25\pi^2\bar\delta_0^3}{81(cL)^2}$.
\end{proof}

\subsection{Global work bound}
\label{app:proof_exact_work}

\begin{lemma}[Work bound]
\label{lem:work}
Under (A1)--(A5) with exact published curvature, $\rho_k = cLw_k$ and
\begin{equation}
    c \;\ge\; c_H(\gamma) = \tfrac53(\Cov + 9),
    \label{eq:c_condition}
\end{equation}
one has, for every $T \ge 1$, with $A = A(c,\Cov) = \tfrac{3c}{20} - \tfrac32 - \tfrac{\Cov}{6} \ge \tfrac{c}{20}$,
\begin{equation}
L\,A\,S_T \;\le\; \Ftilde_0,
\label{eq:work-bound-final}
\end{equation}
\begin{equation}
\sum_{k<T}\rho_k r_k^3 \;=\; cL\,S_T \;\le\; 20\,\Ftilde_0 .
\label{eq:work_rho}
\end{equation}
\end{lemma}

\begin{proof}
Insert \eqref{eq:error_bound} into \eqref{eq:descent-main} and reserve $\tfrac{cL}{10}w_k r_k^3$ at every $k$. For $k < h$, \cref{lem:init} with $n = \min(h,T)$ bounds the Lipschitz terms by $\tfrac{3L}{2}\sum w_kr_k^3$ and the constant error by $\tfrac{25\pi^2\bar\delta_0^3}{81(cL)^2}$. For $k \ge h$, \eqref{eq:young_tail} bounds the current-step coefficient by $\tfrac{L}{3}w_k$ and leaves the history $\tfrac{L}{6}\sum_{k\ge h}\sum_{j=d_k}^{k-1}r_j^3 \le \tfrac{L\Cov}{6}\sum_{j\ge a_1}w_jr_j^3$ by \cref{lem:overlap_global}. Summing over $k < T$ and telescoping $f$ with $f(x_T) \ge \fstar$, the positive smoothness terms total at most $L\big(\tfrac32 + \tfrac{\Cov}{6}\big)S_T$, while the descent provides $\big(\tfrac{c}{4} - \tfrac{c}{10}\big)L\,S_T = \tfrac{3c}{20}L\,S_T$. Hence $L\,A\,S_T \le \Ftilde_0$. The condition $A \ge c/20$ is $\tfrac{c}{10} \ge \tfrac32 + \tfrac{\Cov}{6}$, i.e.\ \eqref{eq:c_condition}, and then $cLS_T \le 20\Ftilde_0$.
\end{proof}

The factor $\pi^2/6$ in $\Ftilde_0$ may be replaced by the finite sum $\sum_{k<\min(h,T)}(k+1)^{-2}$, and $\bar\delta_0$ by the actual error at $x_0$. The stated form is deterministic for given initial data and uniform in the startup duration.

\subsection{First-order rate}
\label{app:proof_exact_grad}

\begin{lemma}[Prefix and delay operator]
\label{lem:prefix_op}
With $K_H = 27 + \Cov$,
\begin{equation}
\sum_{k<T}\frac{P_k^3}{w_k^2} \;\le\; K_H\,S_T .
\label{eq:prefix_op}
\end{equation}
\end{lemma}

\begin{proof}
For $k < h$, $P_k = R_k$, $w_k = k+1$, and the first inequality in \eqref{eq:hardy} bounds these terms by $27\sum_{k<h}w_kr_k^3$. For $k \ge h$, the power mean over the $\hat\tau_k$ summands gives $P_k^3 \le \hat\tau_k^2\sum_{j=d_k}^{k-1}r_j^3 \le w_k^2\sum_{j=d_k}^{k-1}r_j^3$, and \cref{lem:overlap_global} bounds the sum over $k \ge h$ by $\Cov S_T$.
\end{proof}

\begin{lemma}[Per-step gradient bound]
\label{lem:grad-step}
Under (A1), with $a \defeq \tfrac{c+1}{2}$,
\begin{equation}
\norm{\grad(x_{k+1})} \;\le\; aLw_k r_k^2 + Lr_kP_k + \bar\delta_0\,\mathbf{1}_{k<h}\,r_k .
\label{eq:grad-step}
\end{equation}
\end{lemma}

\begin{proof}
By Taylor's theorem and (A1), $\norm{\grad(x_{k+1}) - \grad(x_k) - \hess(x_k)s_k} \le \tfrac{L}{2}r_k^2$. The first condition in \eqref{eq:cubic_certificates} gives $\grad(x_{k+1}) = (\hess(x_k) - H_k)s_k - \tfrac{\rho_k}{2}r_ks_k + [\text{remainder}]$, so $\norm{\grad(x_{k+1})} \le \norm{\hess(x_k) - H_k}\,r_k + \tfrac{\rho_k + L}{2}r_k^2$. Use \eqref{eq:error_bound} and $\tfrac{\rho_k+L}{2} \le aLw_k$.
\end{proof}

\begin{proof}[Proof of \cref{thm:exact}, first-order part]
Raise \eqref{eq:grad-step} to the power $3/2$ with $(x+y+z)^{3/2} \le \sqrt3(x^{3/2}+y^{3/2}+z^{3/2})$:
\begin{equation}
\norm{\grad(x_{k+1})}^{3/2} \;\le\; \sqrt3\Big[a^{3/2}L^{3/2}w_k^{3/2}r_k^3 + L^{3/2}r_k^{3/2}P_k^{3/2} + \bar\delta_0^{3/2}\mathbf{1}_{k<h}\,r_k^{3/2}\Big].
\label{eq:grad-32}
\end{equation}
Multiply by $w_k^{-1/2}$ and sum over $k < T$. The first term gives $a^{3/2}L^{3/2}S_T$. By Cauchy--Schwarz and \cref{lem:prefix_op}, the second gives
\[
L^{3/2}\sum_k\big(w_k^{1/2}r_k^{3/2}\big)\Big(\frac{P_k^{3/2}}{w_k}\Big) \;\le\; L^{3/2}\Big(\sum_kw_kr_k^3\Big)^{1/2}\Big(\sum_k\frac{P_k^3}{w_k^2}\Big)^{1/2} \;\le\; \sqrt{K_H}\,L^{3/2}S_T .
\]
For the third, with $n_0 = \min(h,T)$, Cauchy--Schwarz gives
\[
\bar\delta_0^{3/2}\sum_{k<n_0}\big(w_k^{1/2}r_k^{3/2}\big)(k+1)^{-1} \;\le\; \Big(S_T\,\bar\delta_0^{\,3}\sum_{k<n_0}(k+1)^{-2}\Big)^{1/2} \;\le\; c\sqrt{\frac{27}{50A}}\,L^{1/2}\Ftilde_0 ,
\]
using $\bar\delta_0^3\pi^2/6 \le \tfrac{27}{50}(cL)^2\Ftilde_0$, which follows from \eqref{eq:Ftilde}, and $S_T \le \Ftilde_0/(LA)$. With $Z_T = \sum_{k<T}w_k^{-1/2}$ and \eqref{eq:work-bound-final},
\[
Z_T\,\E\!\left[\norm{\grad(x_{\mathrm{out}})}^{3/2}\right] \;\le\; \sqrt3\Big[\frac{a^{3/2} + \sqrt{K_H}}{A} + c\sqrt{\frac{27}{50A}}\,\Big]L^{1/2}\Ftilde_0 = C_0\,L^{1/2}\Ftilde_0 .
\]
Finally $u \mapsto (1+u)^{-1/2}$ is convex, so Jensen's inequality gives $Z_T \ge T(1+\taubar)^{-1/2}$.
\end{proof}

\subsection{Second-order rate}
\label{app:proof_exact_hess}

\begin{proof}[Proof of \cref{thm:exact}, second-order part]
The argument treats the moving-window history sum as an operator on weighted $\ell^3$, rather than replacing each window by its maximum and then averaging, which pays for the delay twice. Let $T > h$, $N = T - h$, and $D_k = r_k + P_k$. For $k \ge h$ the matrix is an exact delayed Hessian, and \eqref{eq:cubic_certificates}, Weyl's inequality and (A1) give $\lmin(\hess(x_{k+1})) \ge \lmin(H_k) - Lr_k - LP_k \ge -\tfrac{\rho_k}{2}r_k - LD_k$, i.e.
\begin{equation}
-\lmin(\hess(x_{k+1})) \;\le\; L\Big(\tfrac{c}{2}w_kr_k + D_k\Big), \qquad k \ge h.
\label{eq:hess_pointwise}
\end{equation}
Since $w_k \le W_+$ for $k \ge h$,
\begin{equation}
\sum_{k=h}^{T-1}(w_kr_k)^3 \;\le\; W_+^2 S_T .
\label{eq:l3_first}
\end{equation}
$D_k$ has at most $w_k \le W_+$ summands, so the power mean gives $D_k^3 \le W_+^2\big(r_k^3 + \sum_{j=d_k}^{k-1}r_j^3\big)$. Each index $j \ge a_1$ appears at most once as a current step and, by \cref{lem:overlap_global}, at most $\Cov\hat\tau_j$ times in a window; indices $j < a_1$ never occur. Since $1 + \Cov\hat\tau_j \le \Cov w_j$,
\begin{equation}
\sum_{k=h}^{T-1}D_k^3 \;\le\; W_+^2\,\Cov\,S_T .
\label{eq:l3_second}
\end{equation}
By Minkowski's inequality in $\ell^3$, the $N$ values $\tfrac{c}{2}w_kr_k + D_k$, $h \le k < T$, have cubic mean at most $W_+^{2/3}S_T^{1/3}N^{-1/3}\big(\tfrac{c}{2} + \Cov^{1/3}\big)$, so some $k^\star \ge h$ attains at most this value. With \eqref{eq:hess_pointwise} and $S_T \le \Ftilde_0/(LA)$ this is the stated bound with $C_2 = A^{-1/3}(\tfrac{c}{2} + \Cov^{1/3})$.
\end{proof}

\paragraph{A bound for every horizon.} Without restricting the index, \eqref{eq:error_bound} gives $-\lmin(\hess(x_{k+1})) \le L\big((\tfrac{c}{2}+1)w_kr_k + P_k\big) + \bar\delta_0\mathbf{1}_{k<h}$ for all $k$. With $W = \max_{k<T}w_k$, the bounds $\sum_k(w_kr_k)^3 \le W^2S_T$ and $\sum_kP_k^3 \le W^2\sum_kP_k^3/w_k^2 \le W^2K_HS_T$ and Minkowski give some $k < T$ with
\begin{equation}
\lmin(\hess(x_{k+1})) \;\ge\; -A^{-1/3}\Big(\tfrac{c}{2} + 1 + K_H^{1/3}\Big)\frac{L^{2/3}\Ftilde_0^{1/3}W^{2/3}}{T^{1/3}} - \bar\delta_0\Big(\frac{\min(h,T)}{T}\Big)^{1/3}.
\label{eq:hess_all_horizon}
\end{equation}
The post-startup statement of \cref{thm:exact} is cleaner and sharper.

\begin{remark}[On the constant $C_2$ and the choice of $c$]
\label{rem:C2}
The displayed coefficient grows like $c^{2/3}$. Some dependence on $c$ is unavoidable, as the following example shows: for $f(x) = -\tfrac12x^2 + \tfrac16|x|^3$ at $x_0 = 0$ with $H_0 = -1$ (so $\delta_0 = 0$) and unit job durations, the first global cubic step has magnitude $2/(cL)$, so for any fixed horizon a large enough $c$ keeps every iterate arbitrarily close to the strict saddle and every Hessian on the path close to $-1$. At $c = c_H$, where $A = c/20 = (\Cov+9)/12$,
\begin{equation}
C_2(c_H,\Cov) = \Big(\tfrac{12}{\Cov+9}\Big)^{1/3}\Big(\tfrac56(\Cov+9) + \Cov^{1/3}\Big) = \mathcal{O}\big(\Cov^{2/3}\big),
\label{eq:C2min}
\end{equation}
which equals $10.73$ for $\Cov = 2$ and $10.96$ for $\Cov = 2.31$.
\end{remark}

\begin{remark}[On the constants]
The multiplier $c$ trades the two constants against each other. At $\Cov = 2$, $C_0$ is minimized near $c \approx 36$ (value $65.5$, with $C_2 \approx 12.6$), while $C_2$ is minimized near $c \approx 19.6$ (value $10.7$, with $C_0 \approx 83.9$). The constants are explicit sufficient bounds, not optimized ones. The first-order statement concerns the sampled output and the second-order statement an index selected by the analysis; they are separate conclusions, not a joint second-order stationarity guarantee at one point.
\end{remark}

\subsection{A single matrix used throughout}
\label{app:single_matrix}

The Hardy argument also covers the case in which no curvature update ever arrives.

\begin{corollary}[Single-matrix stationarity]
\label{cor:single}
Assume (A1)--(A2). Let $H_k = H_0$ for every $k$, with $H_0$ symmetric and $\norm{H_0 - \hess(x_0)} \le \delta_0$. Use exact gradients, $\rho_k = cL(k+1)$ and $c > 10$; no curvature jobs, (A3) or bounded delay are needed. With $B = \tfrac{3c}{20} - \tfrac32$ and $F_{\mathrm{init}} = F_0 + \tfrac{25\pi^2\delta_0^3}{81(cL)^2}$, one has $LB\sum_{k<T}(k+1)r_k^3 \le F_{\mathrm{init}}$ for every $T$, and with output weights $(k+1)^{-1/2}$,
\[
\E\!\left[\norm{\grad(x_{\mathrm{out}})}^{3/2}\right] \;\le\; \sqrt3\Big[\frac{(\tfrac{c+1}{2})^{3/2} + \sqrt{27}}{B} + c\sqrt{\frac{27}{50B}}\,\Big]\frac{L^{1/2}F_{\mathrm{init}}\sqrt{(T+1)/2}}{T}.
\]
\end{corollary}

\begin{proof}
Apply \cref{lem:init} over the whole horizon: there are no tail windows and no overlap terms. The work and first-order arguments of \cref{app:proof_exact_work,app:proof_exact_grad} then hold with the $\Cov$ contributions removed, $K_H = 27$, $A = B$, and $1 + \taubar = (T+1)/2$.
\end{proof}

The stationarity measure is thus $\mathcal{O}(T^{-1/2})$, and $\E\norm{\grad(x_{\mathrm{out}})} = \mathcal{O}(T^{-1/3})$ by Jensen. This is a statement about the randomized output, not convergence of the iterates to a point, and it relies on increasing regularization. It does not guarantee escape from a saddle with a permanently inaccurate positive semidefinite $H_0$: a zero gradient may yield a zero step forever.

% ============================================================================

\section{Proof of Theorem \ref{thm:unified} (Structured Inexactness)}
\label{app:proof_unified}
% ============================================================================

Write $r_k = \norm{s_k}$, $w_k$, $P_k$, $S_T$ as in \cref{app:proof_exact}, and let $u_k = 0$ for $k < h$ and $u_k = r_{a_{i-1}-1}$ for $a_i \le k < a_{i+1}$, $i \ge 2$. By (A1), (A6) and the assumption on the two matrices based at $x_0$,
\begin{equation}
\norm{\hess(x_k) - \Ht_k} \;\le\; LP_k + \alpha u_k + \bar\delta_0\mathbf{1}_{k<h},
\qquad
\norm{\gt_k - \grad(x_k)} \;\le\; \beta r_{k-1}^2 ,
\label{eq:inexact_errors}
\end{equation}
with $r_{-1} = 0$. At most one of $u_k$ and $\mathbf{1}_{k<h}$ is nonzero.

\begin{lemma}[Reuse and causal-gradient charging]
\label{lem:reuse}
For every horizon,
\begin{equation}
\sum_{k<T}u_k^3 \;\le\; \gamma^2 S_T ,
\label{eq:reuse_global}
\end{equation}
and $\sum_{k<T}r_{k-1}^3 \le \sum_{k<T}r_k^3$.
\end{lemma}

\begin{proof}
During job $i \ge 2$ the reference index $\ell = a_{i-1} - 1$ is fixed and is used $\Delta_i$ times. For $i = 2$, $w_\ell = \Delta_0$ and $\Delta_2 \le \gamma^2\Delta_0$. For $i \ge 3$, $\ell$ lies in job $i-2$ and $\tau_\ell = \Delta_{i-3} + \Delta_{i-2} - 1 \ge \Delta_{i-2}$, because durations are positive integers, so $\Delta_i \le \gamma^2\Delta_{i-2} \le \gamma^2 w_\ell$. The reference indices are distinct, which gives \eqref{eq:reuse_global}. The second claim is a shift of a nonnegative finite sum with $r_{-1} = 0$.
\end{proof}

\subsection{Work bound}

\begin{lemma}[Work bound, structured inexactness]
\label{lem:work-inexact}
Under the assumptions of \cref{thm:unified},
\begin{equation}
E_T \;\defeq\; \sum_{k<T}\rho_k r_k^3 \;\le\; 20\,\Ftilde_0 .
\label{eq:work-inexact-final}
\end{equation}
\end{lemma}

\begin{proof}
First-order optimality for the sampled model, $\gt_k + \Ht_ks_k + \tfrac{\rho_k}{2}r_ks_k = 0$, and the global-minimizer condition $\Ht_k + \tfrac{\rho_k}{2}r_kI \succeq 0$ give, as in \cref{lem:descent},
\[
f(x_{k+1}) - f(x_k) \le \tfrac{1}{2}\inner{(\hess(x_k) - \Ht_k) s_k}{s_k} + \inner{\grad(x_k) - \gt_k}{s_k} + \Big(\tfrac{L}{6} - \tfrac{\rho_k}{4}\Big)r_k^3.
\]
Insert \eqref{eq:inexact_errors}. The Lipschitz part and the constant error $\bar\delta_0$ are handled exactly as in \cref{lem:work}, reserving only $\tfrac{cL}{10}w_k$ for the latter; they contribute $L(\tfrac32 + \tfrac{\Cov}{6})S_T$ and the additive $\tfrac{25\pi^2\bar\delta_0^3}{81(cL)^2}$. For the other errors,
\[
\tfrac{\alpha}{2}u_kr_k^2 \le \tfrac{\alpha}{6}u_k^3 + \tfrac{\alpha}{3}r_k^3,
\qquad
\beta r_{k-1}^2r_k \le \tfrac{2\beta}{3}r_{k-1}^3 + \tfrac{\beta}{3}r_k^3 .
\]
Summing and using \cref{lem:reuse}, the gradient error costs at most $\beta\sum_kr_k^3$ in total and the reuse term at most $\tfrac{\gamma^2\alpha}{6}\sum_kw_kr_k^3$. Since $\tfrac{\rho_k}{4} = \tfrac{c}{4}Lw_k + (1 + \tfrac{\gamma^2}{4})\alpha + \tfrac{\gamma^2\alpha\hat\tau_k}{2} + \tfrac{3\beta}{2}$, the residual coefficient of $r_k^3$ after telescoping is at least
\[
R_k \;\defeq\; L\,A\,w_k + \Big(\tfrac23 + \tfrac{\gamma^2}{12}\Big)\alpha + \tfrac{\gamma^2\alpha\hat\tau_k}{3} + \tfrac{\beta}{2}:
\]
the $\alpha$ part before reuse is $(1 + \tfrac{\gamma^2}{4})\alpha + \tfrac{\gamma^2\alpha\hat\tau_k}{2} - \tfrac{\alpha}{3}$, and subtracting $\tfrac{\gamma^2\alpha w_k}{6}$ gives the displayed expression. Hence $\sum_kR_kr_k^3 \le \Ftilde_0$. Since $A \ge c/20$, the ratios of the coefficients of $\rho_k$ to those of $R_k$ are $20$, $\tfrac{12(4+\gamma^2)}{8+\gamma^2} \le 12$, $6$ and $12$, so $\rho_k \le 20R_k$ and $E_T \le 20\Ftilde_0$.
\end{proof}

\subsection{First-order rate for Theorem \ref{thm:unified}}

\begin{proof}[Proof of \cref{thm:unified}, first-order part]
A Taylor expansion of the gradient and first-order optimality give
\[
\norm{\grad(x_{k+1})} \;\le\; \tfrac{\rho_k+L}{2}r_k^2 + Lr_kP_k + (\alpha u_k + \bar\delta_0\mathbf{1}_{k<h})r_k + \beta r_{k-1}^2 ,
\]
and $\tfrac{\rho_k + L}{2} \le \rho_k$. Group the right-hand side into four nonnegative terms (leading, history, curvature error, gradient error), use $(z_1 + \dots + z_4)^{3/2} \le 2\sum_iz_i^{3/2}$, multiply by $\rho_k^{-1/2}$ and sum over $k < T$.
\begin{itemize}[itemsep=0.1em]
    \item The leading sum is $E_T$.
    \item Since $\rho_k \ge cLw_k$, Cauchy--Schwarz and \cref{lem:prefix_op} bound the history sum by $L\sqrt{K_H}\,S_T/\sqrt{c} \le \sqrt{K_H}\,E_T/c^{3/2}$.
    \item Since $\rho_k \ge 4\alpha$, the structured part is at most $\tfrac{\alpha}{4}\sum_k(u_k^3 + r_k^3)$; using $\rho_k \ge 2\alpha w_k$, $\rho_k \ge 4\alpha$ and \eqref{eq:reuse_global}, this is at most $(\tfrac{\gamma^2}{8} + \tfrac1{16})E_T$.
    \item The constant initial error contributes at most $\big(E_T\sum_{k<h}\bar\delta_0^3/\rho_k^2\big)^{1/2} \le \sqrt{54/5}\,\Ftilde_0$ by Cauchy--Schwarz, $\rho_k \ge cL(k+1)$ and \eqref{eq:Ftilde}.
    \item Since $\rho_k \ge 6\beta$, the gradient-error sum is at most $\sum_k\beta^{3/2}r_{k-1}^3/\sqrt{\rho_k} \le \tfrac{\beta}{\sqrt6}\sum_kr_k^3 \le \tfrac{E_T}{6\sqrt6}$.
\end{itemize}
For a zero parameter the corresponding term is omitted. With $E_T \le 20\Ftilde_0$ and $Z_T = \sum_k\rho_k^{-1/2} \ge T\bar\rho^{-1/2}$ by Jensen, since $\rho_k$ is affine in $\hat\tau_k$, this gives the stated bound with $C_1 = 2\big[20(\tfrac{17}{16} + \tfrac{\gamma^2}{8} + \tfrac{\sqrt{K_H}}{c^{3/2}} + \tfrac1{6\sqrt6}) + \sqrt{54/5}\big]$.
\end{proof}

\subsection{Second-order rate for Theorem \ref{thm:unified}}

\begin{proof}[Proof of \cref{thm:unified}, second-order part]
For $k \ge h$, \eqref{eq:cubic_certificates} with $\Ht_k$, Weyl's inequality and \eqref{eq:inexact_errors} give $-\lmin(\hess(x_{k+1})) \le \tfrac{\rho_k}{2}r_k + LD_k + \alpha u_k$ with $D_k = r_k + P_k$; the gradient error does not enter the eigenvalue bound, and no initial-error term appears in this index set. Minkowski in $\ell^3$ over the $N$ indices $h \le k < T$, with $\sum_k(\rho_kr_k)^3 \le \rho_+^2E_T$, $\sum_kD_k^3 \le W_+^2\Cov S_T$ from \eqref{eq:l3_second}, $\sum_ku_k^3 \le \gamma^2S_T$, and $S_T \le E_T/(cL) \le 20\Ftilde_0/(cL)$, gives the stated bound.
\end{proof}

% ============================================================================
% ============================================================================
\section{Structured Quasi-Newton curvature}
\label{app:qn_full}
% ============================================================================

This appendix records what exact secant equations imply for the curvature error. It is not used elsewhere in the paper and gives no convergence rate for a Quasi-Newton oracle.

\subsection{Secant and subspace estimates}

Prior treatments of Quasi-Newton within cubic regularization \citep{kamzolov2023cubic} work with a uniform error bound $\norm{\HBFGS_k - \hess(x_k)} \le \delta$, compensated by a quadratic correction to the cubic model, in the convex and star-convex settings. We consider a symmetric matrix $\HBFGS_k$ satisfying \emph{exact} secant equations $\HBFGS_k s_i = y_i$ for $i$ in a memory $\mathcal{M}_k \subseteq \{k-m, \dots, k-1\}$ (truncated at $0$), with $y_i = \grad(x_{i+1}) - \grad(x_i)$. Exact multi-secant equations together with a spectral bound are assumptions on the oracle. A multi-secant update \citep{schnabel1983} produces such a matrix when the block curvature matrix $S_k^\top Y_k$ is symmetric positive definite, but compatibility with fixed bounds $\mu I \preceq \HBFGS_k \preceq MI$ additionally requires, for instance, $\mu S_k^\top S_k \preceq S_k^\top Y_k \preceq M S_k^\top S_k$, and we do not claim that a given construction always meets these safeguards. Standard L-BFGS in general preserves only its latest secant condition, not the whole memory; for it, only the single-direction statement of \cref{lem:secant} with $i = k-1$ applies.

\begin{lemma}[Secant Error]
\label{lem:secant}
Under (A1), if $\HBFGS_k s_i = y_i$ for some $i \le k-1$, then
\begin{equation}
    \norm{(\HBFGS_k - \hess(x_k))\, s_i} \;\le\; L\,\norm{s_i}\,\sum_{j=i}^{k-1} \norm{s_j}.
\end{equation}
\end{lemma}

\begin{proof}
Since $y_i = \int_0^1\hess(x_i + ts_i)s_i\,dt$, the secant equation gives $(\HBFGS_k - \hess(x_k))s_i = \int_0^1[\hess(x_i+ts_i) - \hess(x_k)]s_i\,dt$, and $\norm{x_k - (x_i + ts_i)} \le \sum_{j=i}^{k-1}\norm{s_j}$ for $t \in [0,1]$.
\end{proof}

If damping is used, so that $\HBFGS_k s_i = y_i + d_i$, the bound of \cref{lem:secant} gains $\norm{d_i}$, and the subspace estimate below gains a corresponding residual on the secant subspace itself; it cannot be absorbed into the off-subspace error. We state the results for exact secants.

\begin{proposition}[Subspace and operator estimates]
\label{prop:lbfgs_fits}
Assume (A1), exact secant equations on $\mathcal{M}_k$, and:
\begin{enumerate}[label=(QN\arabic*),itemsep=0em]
    \item \textbf{Bounded curvature at the iterates:} $\norm{\hess(x_k)} \le H$ at every iterate actually visited, and $\norm{\HBFGS_k} \le M$;
    \item \textbf{Secant-basis conditioning:} $S_k = [s_i]_{i\in\mathcal{M}_k}$ has full column rank, and $\kappa_S \defeq s_M/\sigma_k \le \bar\kappa_S < \infty$, where $s_M$ is the largest of the norms $\norm{s_i}$ and $\sigma_k > 0$ the smallest singular value of $S_k$.
\end{enumerate}
With $R_k = \sum_{j=k-m}^{k-1}\norm{s_j}$ (lower limit truncated at $0$) and $V_k = \mathrm{span}\{s_i : i\in\mathcal{M}_k\}$, for every $v \in V_k$,
$\norm{(\HBFGS_k - \hess(x_k))v} \le LR_k\sqrt{m}\,\kappa_S\norm{v}$, and $\norm{\HBFGS_k - \hess(x_k)} \le M + H \eqqcolon \delta_u$.
\end{proposition}

\begin{proof}
Write $v = S_ka$ and let $D_k$ be diagonal with entries $\norm{s_i}$. By \cref{lem:secant} and the triangle inequality, $\norm{(\HBFGS_k - \hess(x_k))v} \le LR_k\norm{D_ka}_1$. Cauchy--Schwarz over the $m' \le m$ coordinates gives $\norm{D_ka}_1 \le \sqrt m\,\norm{D_ka}_2 \le \sqrt m\,s_M\norm{a}_2$, and $\sigma_k\norm{a}_2 \le \norm{S_ka}_2 = \norm{v}$. The operator bound follows from (QN1) and the triangle inequality.
\end{proof}

\begin{remark}[Why $\kappa_S$ is assumed rather than derived]
\label{rem:kappa_honest}
The conditioning of the secant basis is a separate quantity from the spectral spread $M/\mu$: nearly parallel memory directions make $\sigma_k$ arbitrarily small whatever the spectrum of $\HBFGS_k$. It is also not equivalent to the curvature test on $s_i^\top y_i$: for $f(x) = \tfrac12\norm{x}^2$ with $s_1 = e_1$, $s_2 = e_1 + \eta e_2$ and $y_i = s_i$, both curvature ratios equal $1$ while $\sigma_k \to 0$ as $\eta \to 0$. Controlling $\kappa_S$ requires a direction-conditioning safeguard, for instance a Gram or QR test on the memory.
\end{remark}

\subsection{What the estimates do and do not give}

Decomposing $s_k$ along $V_k$ and its complement, with angle $\theta_k$ between $s_k$ and $V_k$, gives
\begin{equation}
    \norm{(\HBFGS_k - \hess(x_k))\, s_k} \;\le\; \big(L R_k \sqrt m\,\kappa_S\, \cos\theta_k + \delta_u \sin\theta_k\big)\, \norm{s_k}
    \label{eq:hybrid_per_step}
\end{equation}
at the point and for the matrix to which \cref{prop:lbfgs_fits} applies. This is a directional, history-dependent estimate, not the operator-norm form of (A6). In the asynchronous method the matrix in use at iteration $k$ was formed from a memory at an older snapshot and is reused for a whole job, which requires its own error accounting, indexed by jobs. Applying the uniform-error analysis with $\delta = \delta_u$, which requires this bound at each dispatch snapshot and an initial matrix $H_0$ satisfying it, gives the bound of \cref{thm:robust}; it does not exploit the angles $\theta_k$. An empty memory has secant subspace $\{0\}$, and then only the uniform bound applies. This upper bound does not establish a lower bound excluding faster convergence under additional directional structure. A convergence analysis exploiting the directional estimate \eqref{eq:hybrid_per_step} for asynchronously reused matrices remains open.

% ============================================================================
\section{Sampling: probability and budget}
\label{app:budget_proof}
% ============================================================================

\subsection{Realizing (A6) with high probability}
\label{app:a6_realize}

With the counts \eqref{eq:batch_sizes}, every tolerance is measurable before its batch is drawn, so matrix Hoeffding applies conditionally on the past and each call satisfies its bound in (A6) with probability at least $1 - p$. Calls that evaluate the full sum, and reuses of an already published matrix, create no new failure event.

\paragraph{Finite horizon.} With at most $T$ gradient calls and $T$ curvature jobs, choose $p = \eta/(2T)$ for every call. A conditional union bound gives a good event $\mathcal{G}$, on which (A6) holds throughout the run, with $\Pr(\mathcal{G}) \ge 1 - \eta$ for $0 < \eta < 1$; for $T > 2$, $\eta = 2/T$ recovers $p = T^{-2}$.

\paragraph{Anytime.} Numbering the non-exact calls $j = 0, 1, \dots$ and using $p_j = 6\eta/(\pi^2(j+1)^2)$, the union bound covers every finite prefix and the whole infinite run.

\paragraph{From the good event to an expectation.} On $\mathcal{G}$, \cref{thm:unified} gives $\E_I[\norm{\grad(x_{\mathrm{out}})}^{3/2} \mid \text{run}] \le C_1\Ftilde_0\sqrt{\bar\rho(\text{run})}/T$, where $\E_I$ averages over the output index only. If in addition $\norm{\grad(x)} \le G$ on a region containing every possible trajectory, including those on which an oracle fails, then
\begin{equation}
\E\!\left[\norm{\grad(x_{\mathrm{out}})}^{3/2}\right] \;\le\; \frac{C_1\Ftilde_0}{T}\,\E\big[\mathbf{1}_{\mathcal{G}}\sqrt{\bar\rho}\,\big] + \eta\,G^{3/2}.
\label{eq:unconditional}
\end{equation}
A deterministic bound on $\bar\rho$ may replace the expectation, but a random realized average age cannot simply be left outside an unconditional expectation. A bounded sublevel set is not automatically invariant here, since the descent estimate carries positive history terms and no acceptance rule is used; without a valid bound on failed trajectories, only the high-probability conditional statement on $\mathcal{G}$ holds.

\subsection{The aggregate budget}

\begin{proof}[Proof of \cref{prop:budget}]
Each call evaluates at most $n$ components by \eqref{eq:batch_sizes}. There is one gradient call per iteration, and curvature jobs are dispatched at distinct indices $a_i < T$ because durations are positive integers. Summing gives \eqref{eq:budget}.
\end{proof}

% ============================================================================
\section{Proof of Theorem \ref{thm:robust} (Robustness)}
\label{app:proof_robust}
% ============================================================================

For $\delta = 0$ the noise terms vanish and the exact-curvature argument of \cref{app:proof_exact} applies, so assume $\delta > 0$ below. We prove the anytime variant; the fixed-horizon variant is identical with $\sqrt{k+1}$ replaced by $\sqrt{T}$. Write $r_k$, $w_k$, $P_k$, $S_T$ as in \cref{app:proof_exact}. Every matrix has error at most $\delta$ at its base point, so $\norm{\hess(x_k) - \Ht_k} \le \delta + LP_k$ for every $k$.

\subsection{Work bound under uniform noise}

\begin{lemma}[Anytime work bound]
\label{lem:work-anytime}
Under the assumptions of \cref{thm:robust},
\begin{equation}
\sum_{k<T}\rho_k r_k^3 \;\le\; \mathcal{K}(T) \;=\; 20\,\Ftilde_0 + \frac{1000\,(1 + \ln T)}{27\,\nu^2}.
\label{eq:work-anytime}
\end{equation}
\end{lemma}

\begin{proof}
In \eqref{eq:descent-main} reserve $\rho_k/10$ and use $\tfrac{\delta}{2}r_k^2 \le \tfrac{\rho_k}{10}r_k^3 + \tfrac{50\delta^3}{27\rho_k^2}$. The Lipschitz part is handled as in \cref{lem:work}: \cref{lem:init} for $k < h$ and \eqref{eq:young_tail} with \cref{lem:overlap_global} for $k \ge h$ give a total of at most $L(\tfrac32 + \tfrac{\Cov}{6})S_T$. Since $\rho_k \ge cLw_k$ and $c \ge c_H$ give $L(\tfrac32 + \tfrac{\Cov}{6})w_k \le \tfrac{\rho_k}{10}$, the residual coefficient satisfies $\tfrac{3\rho_k}{20} - L(\tfrac32 + \tfrac{\Cov}{6})w_k \ge \tfrac{\rho_k}{20}$. Telescoping $f$ with $f(x_T) \ge \fstar$,
\[
\tfrac1{20}\sum_k \rho_k r_k^3 \;\le\; F_0 + \tfrac{50}{27}\sum_{k<h}\frac{\delta^3}{(cL)^2(k+1)^2} + \tfrac{50}{27}\sum_{k\ge h}\frac{\delta^3}{\nu^2\delta^3(k+1)} \;\le\; \Ftilde_0 + \frac{50\,(1+\ln T)}{27\,\nu^2},
\]
using $\rho_k \ge cL(k+1)$ for $k < h$, $\rho_k^2 \ge \nu^2\delta^3(k+1)$ for $k \ge h$, and the harmonic sum.
\end{proof}

\subsection{Rate}

\begin{proof}[Proof of \cref{thm:robust}]
A Taylor expansion and first-order optimality give $\norm{\grad(x_{k+1})} \le \rho_kr_k^2 + Lr_kP_k + \delta r_k$, using $\tfrac{\rho_k+L}{2} \le \rho_k$. Raise to the power $3/2$ with $(x+y+z)^{3/2} \le \sqrt3(x^{3/2}+y^{3/2}+z^{3/2})$, multiply by $w_k^{-1/2}$ and sum; $Z_T = \sum_kw_k^{-1/2} \ge T/\sqrt{1+\taubar}$ by Jensen.

\emph{Work term.} $\sum_kw_k^{-1/2}\rho_k^{3/2}r_k^3 = \sum_k\rho_kr_k^3(\rho_k/w_k)^{1/2} \le (\sqrt{cL} + \sqrt\nu\,\delta^{3/4}T^{1/4})\mathcal{K}(T)$, since $\rho_k/w_k \le \max(cL, \nu\delta^{3/2}\sqrt{k+1})$.

\emph{History term.} By Cauchy--Schwarz and \cref{lem:prefix_op}, $L^{3/2}\sum_kw_k^{-1/2}(r_kP_k)^{3/2} \le \sqrt{K_H}L^{3/2}S_T \le \tfrac{\sqrt{K_H}}{c}L^{1/2}\mathcal{K}(T)$, since $cLS_T \le \mathcal{K}(T)$.

\emph{Startup noise, $k < h$.} By Cauchy--Schwarz, $\delta^{3/2}\sum_{k<h}w_k^{-1/2}r_k^{3/2} \le (S_T\,\delta^3\pi^2/6)^{1/2}$, and $\delta^3\pi^2/6 \le \tfrac{27}{50}(cL)^2\Ftilde_0 \le \tfrac{27}{1000}(cL)^2\mathcal{K}(T)$, so the term is at most $\sqrt{27/1000}\,\sqrt{cL}\,\mathcal{K}(T) \le \tfrac15\sqrt{cL}\,\mathcal{K}(T)$.

\emph{Tail noise, $k \ge h$.} Dropping $w_k^{-1/2} \le 1$, Cauchy--Schwarz with $\sum_k\rho_kr_k^3$ and $\sum_{k<T}\rho_k^{-1} \le \tfrac{2\sqrt T}{\nu\delta^{3/2}}$ give at most $\delta^{3/4}\sqrt{2\mathcal{K}(T)/\nu}\;T^{1/4}$.

Dividing by $Z_T$ gives the stated bound. For the fixed-horizon schedule $\rho_k = \max(cLw_k, \nu\delta^{3/2}\sqrt T)$ one has $\sum_{k\ge h}\delta^3/\rho_k^2 \le \nu^{-2}$ and $\sum_k\rho_k^{-1} \le \sqrt T/(\nu\delta^{3/2})$, which replaces $1+\ln T$ by $1$ in $\mathcal{K}(T)$.
\end{proof}

% ============================================================================
\section{Elapsed time with arbitrary iteration costs}
\label{app:phys_time}
% ============================================================================

Let $g_k^{\mathrm{V}}$ be the full non-curvature cost of iteration $k$ of Vanilla and $c_k^{\mathrm{V}}$ its curvature-job cost, and $g_k^{\mathrm{A}}$ the full cost of iteration $k$ of Split-Client, including its cubic solve with the published factorization. Exactly,
\[
    W_\mathrm{V} = \sum_{k<T_\mathrm{V}}(g_k^\mathrm{V} + c_k^\mathrm{V}),\qquad W_\mathrm{A} = \sum_{k<T_\mathrm{A}} g_k^\mathrm{A},\qquad
    \frac{W_\mathrm{V}}{W_\mathrm{A}} = \frac{T_\mathrm{V}}{T_\mathrm{A}}\,(1 + r_\mathrm{V})\,\frac{\bar g_\mathrm{V}}{\bar g_\mathrm{A}},\quad r_\mathrm{V} \defeq \frac{\sum_k c_k^\mathrm{V}}{\sum_k g_k^\mathrm{V}}.
\]
Equal iteration costs do not make $r_\mathrm{V}$ equal to the average age $\taubar$ of the Split-Client run: $r_\mathrm{V}$ is a job cost, $\taubar$ an average age, and for constant duration $D$ they are $D$ and $(3D-1)/2$ respectively. Two upper bounds on iteration counts do not give an equality of runtimes, so we compare upper-bound scalings only, as in \cref{sec:wall_clock_comparison}.

A single-run elapsed-time bound nevertheless holds for arbitrary costs $g_k > 0$. Suppose a theorem gives $Z_T\,\E\norm{\grad(x_{\mathrm{out}})}^{3/2} \le B$ with $Z_T = \sum_k\rho_k^{-1/2}$, and let $W = \sum_{k<T}g_k$. Writing $W = \sum_k(\rho_k^{-1/4})(g_k\rho_k^{1/4})$, Cauchy--Schwarz gives $Z_T \ge W^2/\sum_k g_k^2\rho_k^{1/2}$, hence
\[
    \E\!\left[\norm{\grad(x_{\mathrm{out}})}^{3/2}\right] \;\le\; B\,\frac{\sum_k g_k^2\,\rho_k^{1/2}}{W^2}.
\]
For \cref{thm:exact}, $B = C_0\Ftilde_0/\sqrt{c}$; for \cref{thm:unified}, $B = C_1\Ftilde_0$. The weight $g_k^2$ makes the bound invariant under a change of time unit. The conversion does not establish a cross-method speedup ratio.

\paragraph{What an iteration costs.} Let $G$, $S$ and $C$ be the costs of a gradient pass, of the remaining per-iteration work (cubic solve and update), and of a curvature job. Two ratios are relevant: $q_{\mathrm{grad}} = C/G$ and the job duration in iteration units $D = C/(G+S)$, which is the quantity used in \cref{def:timeline} and \cref{sec:wall_clock_comparison}. In the dense model $G = \Theta(nd)$, $C = \Theta(nd^2 + d^3)$ and $S = \Theta(d^2)$ for the spectral transforms, so $q_{\mathrm{grad}} = \Theta(d + d^2/n)$ while $D = \Theta(d)$. Measured ratios can differ from both, depending on hardware, contention and implementation. A reusable eigensystem supports shifted cubic solves; a single unshifted Cholesky factor does not in general solve every shifted indefinite cubic model.

% ============================================================================
\section{Average versus maximum delay}
\label{app:avg_vs_max}
% ============================================================================

\begin{remark}[Average versus maximum delay: an illustration]
\label{rem:heavy_tail}
The first-order bound of \cref{thm:exact} is governed by the \emph{realized} average delay $\taubar = \tfrac1T\sum_k \tau_k$, through the Jensen step $Z_T \ge T/\sqrt{1+\taubar}$ in \cref{app:proof_exact_grad}, and never by $\tmax$. To build intuition for why this matters, consider three schematic delay profiles with the same \emph{nominal} mean $10$: constant, $\tau_k \equiv 10$; two-point, $\tau_k \in \{5,55\}$ in proportion $(0.9,0.1)$; and heavy-tailed, most delays small with a rare few very large (e.g.\ i.i.d.\ Pareto with tail index $2$, finite mean, infinite variance). A $\tmax$-based analysis would penalize the latter two by factors of $\sqrt{55/10} \approx 2.3$ and, since the maximum of $T$ such Pareto draws grows like $\sqrt{T}$, by an unbounded amount; a $\taubar$-based bound is governed by the average.

This is an algebraic comparison of delay profiles, not a robustness result for i.i.d.\ heavy-tailed job durations: (A3) and (A4) do not hold uniformly for an infinite sequence of unbounded Pareto durations, independent Pareto ages need not be realizable as the sawtooth age process of a serial queue, and on a finite run the realized $\gamma$ and $\tmax$ can grow with the outliers. The rigorous content is that the first-order bound of \cref{thm:exact} involves whatever $\taubar$ actually occurs, with constants depending on $\gamma$.
\end{remark}

\section{Experimental Details}
\label{app:additional_experiments}
% ============================================================================

This appendix provides additional details on the experimental setup of \cref{sec:experiments}. As in the main text, ``Async'' denotes Split-Client. All experiments were run on Google Colab.

\subsection{Implementation}
\label{app:implementation}

The Split-Client algorithm is implemented in Python using the \texttt{multiprocessing} module (CPython's GIL would otherwise serialize gradient and curvature work). The curvature client runs as a separate OS process and writes the new factorization into pre-allocated shared-memory buffers (\texttt{multiprocessing.shared\_memory}); a version counter signals the gradient client to load the new matrix at the next iteration. Reads and publications are protected: each job works on a copy of the iterate taken under a lock, and a factorization becomes visible to the gradient client only after it is completely written into a buffer the gradient client is not reading, with the version counter updated last. The matrix used at every iteration is therefore exactly the output of one completed job, as the analysis assumes. Memory profiling uses \texttt{tracemalloc} for the gradient process and \texttt{psutil} (resident set size) for the worker process; \cref{fig:async-memory} reports their sum. These are different metrics, and shared pages may be counted in both, so the factor of two is approximate.

\subsection{Additional Experiments on synthetic data}

\subsubsection{Geman--McClure regression}

The Geman--McClure-regularized least-squares objective:
\[
    f(x) = \tfrac{1}{2n}\norm{Ax - b}^2 + \lambda\sum_{i=1}^d \frac{x_i^2}{1+x_i^2}.
\]
The penalty $\phi(t) = t^2/(1+t^2)$ is smooth, with Lipschitz second derivative, and is a smooth surrogate for the $\ell_0$ ``norm''. The objective contains this nonconvex penalty, but whether the full objective is nonconvex depends on the least-squares curvature and the regularization strength: since $\phi''(t) = (2 - 6t^2)/(1+t^2)^3 \ge -\tfrac12$, with equality at $t = \pm1$, one has $\hess(x) \succeq A^\top A/n - \tfrac{\lambda}{2}I$, and the objective is nonconvex exactly when $\lmin(A^\top A/n) < \lambda/2$, strongly convex when $\lmin(A^\top A/n) > \lambda/2$. We do not classify the generated instances here; in either case they test the computational benefit of the architecture. Since $A^\top A/n$ and $A^\top b/n$ do not depend on $x$, they could be precomputed once, after which a curvature job reduces to a diagonal update and a dense factorization, and a gradient costs $\Theta(d^2)$. Our implementation does not use this caching: every curvature job assembles the Hessian in full, in all Geman--McClure experiments and in \cref{fig:speedup-vs-d,fig:time-profile-wide}. The reported timings therefore include assembly costs that a cached implementation would avoid, and the measured advantage over Vanilla would likely be smaller with caching.

\paragraph{Data and parameters.} All results are averaged over 5 random seeds. $A \in \R^{n \times d}$ with i.i.d.\ $\mathcal{N}(0, 1)$ entries, $n = 5000$, $d = 1000$. Ground truth $x_{\mathrm{true}}$ has 10\% non-zero entries drawn from $\mathcal{N}(0, 1)$; targets $b = A x_{\mathrm{true}} + \epsilon$ with $\epsilon \sim \mathcal{N}(0, 10^{-2} I)$. Regularization $\lambda = 10^{-2}$. Initialization: $x_0 = 0$.

\paragraph{Hyperparameters used in \cref{fig:nonconvex-loss-appendix}.}
Async: $\rho = 1.0\times 10^4$. Lazy: $\rho = 4.0\times 10^3$, $m = 100$. Vanilla: $\rho = 4.0\times 10^3$. Each $\rho$ was tuned by a search over powers of $10$ followed by a finer search around the best value, and is the best among the tested values. In the other experiments Lazy uses the reuse period $m = d$, the choice suggested by the arithmetic-cost analysis of \citet{doikov2023lazy}; on this problem $m = 100$ performed better than $m = d = 1000$, so we report it.

\subsubsection{One-layer tanh regression}

A nonconvex objective serving as a proxy for single-layer neural network training:
\[
    f(x) = \tfrac{1}{2n}\sum_{i=1}^n (\tanh(a_i^\top x) - y_i)^2.
\]
The gradient is $\nabla f(x) = \tfrac{1}{n} A^\top \!\big[(\tanh(Ax) - y) \odot (1 - \tanh^2(Ax))\big]$. The Hessian $\nabla^2 f(x) = \tfrac{1}{n} A^\top \mathrm{diag}(\alpha) A$ has diagonal weights
\[
    \alpha_i = (1-\tanh^2(z_i))^2 \;-\; 2\tanh(z_i)\,(\tanh(z_i) - y_i)\,(1-\tanh^2(z_i)), \qquad z_i = a_i^\top x,
\]
which are not all positive: when the residual term dominates, $\alpha_i$ can be negative, in which case $\nabla^2 f$ can be indefinite.

\paragraph{Subproblem solver.} For each cubic subproblem, we use a rank-$r$ approximation of the Hessian computed by randomized SVD, with $r = \min(50, d)$. We do not certify that this approximation retains the signs of all eigenvalues, how neglected directions enter the model, or that the resulting cubic subproblem, including its hard case, is solved globally; these experiments therefore evaluate an implementation variant rather than the exact algorithm of the theorems.

\paragraph{Data and parameters.} All results are averaged over 5 random seeds. $A$ Gaussian with $n = 1000$, $d = 500$. Ground truth $x_{\mathrm{true}}$ standard Gaussian; $y_i = \tanh(a_i^\top x_{\mathrm{true}}) + \epsilon$ with $\epsilon \sim \mathcal{N}(0, 10^{-3})$.

\begin{figure}[h]
  \centering
  \begin{minipage}[b]{.49\linewidth}
    \centering
    \includegraphics[width=\linewidth]{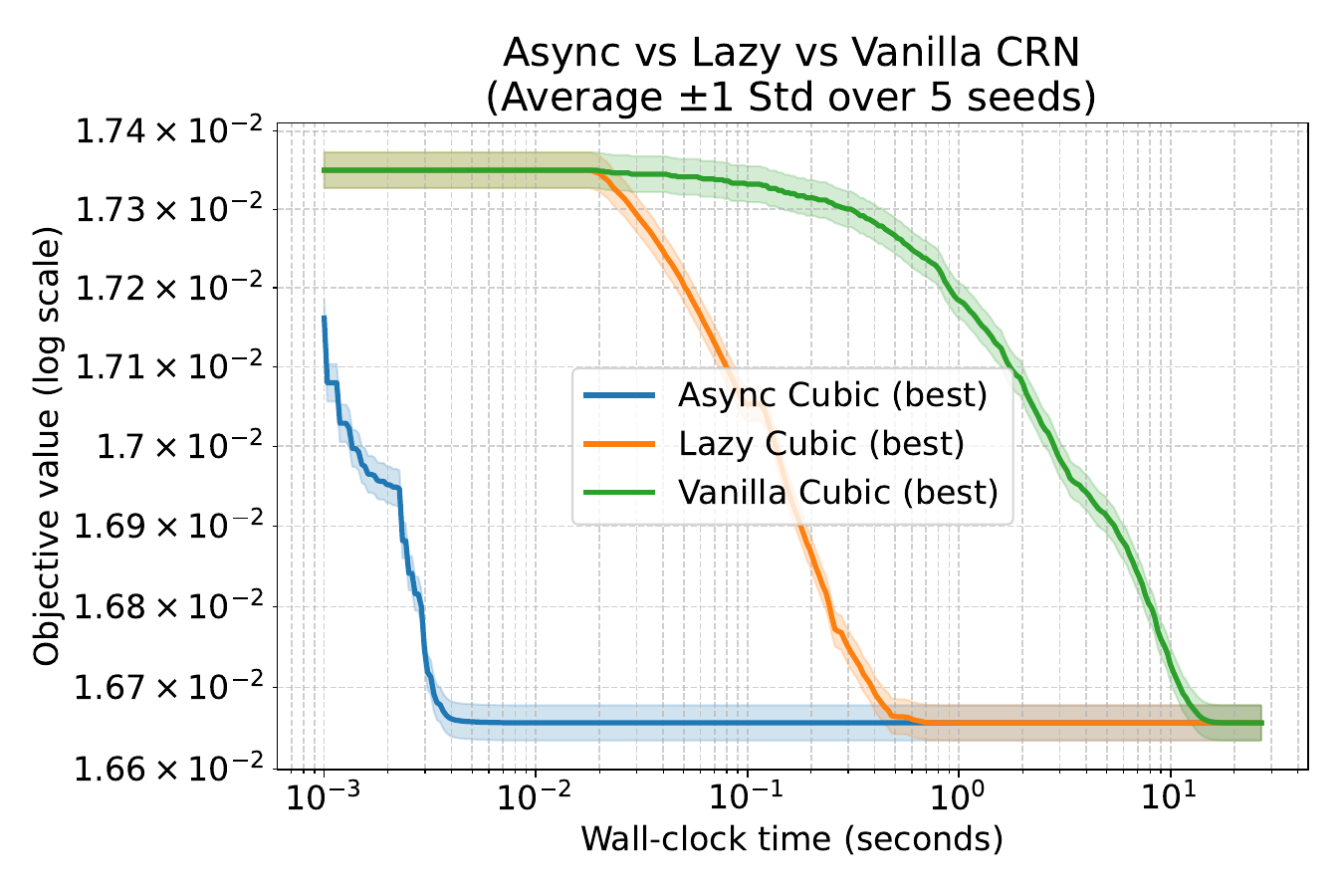}
    \caption{Geman--McClure regression ($d{=}1000$): objective versus wall-clock time, averaged over 5 random seeds. The shaded regions denote $\pm 1$ standard deviation.}
    \label{fig:nonconvex-loss-appendix}
  \end{minipage}\hfill
  \begin{minipage}[b]{.49\linewidth}
    \centering
    \includegraphics[width=\linewidth]{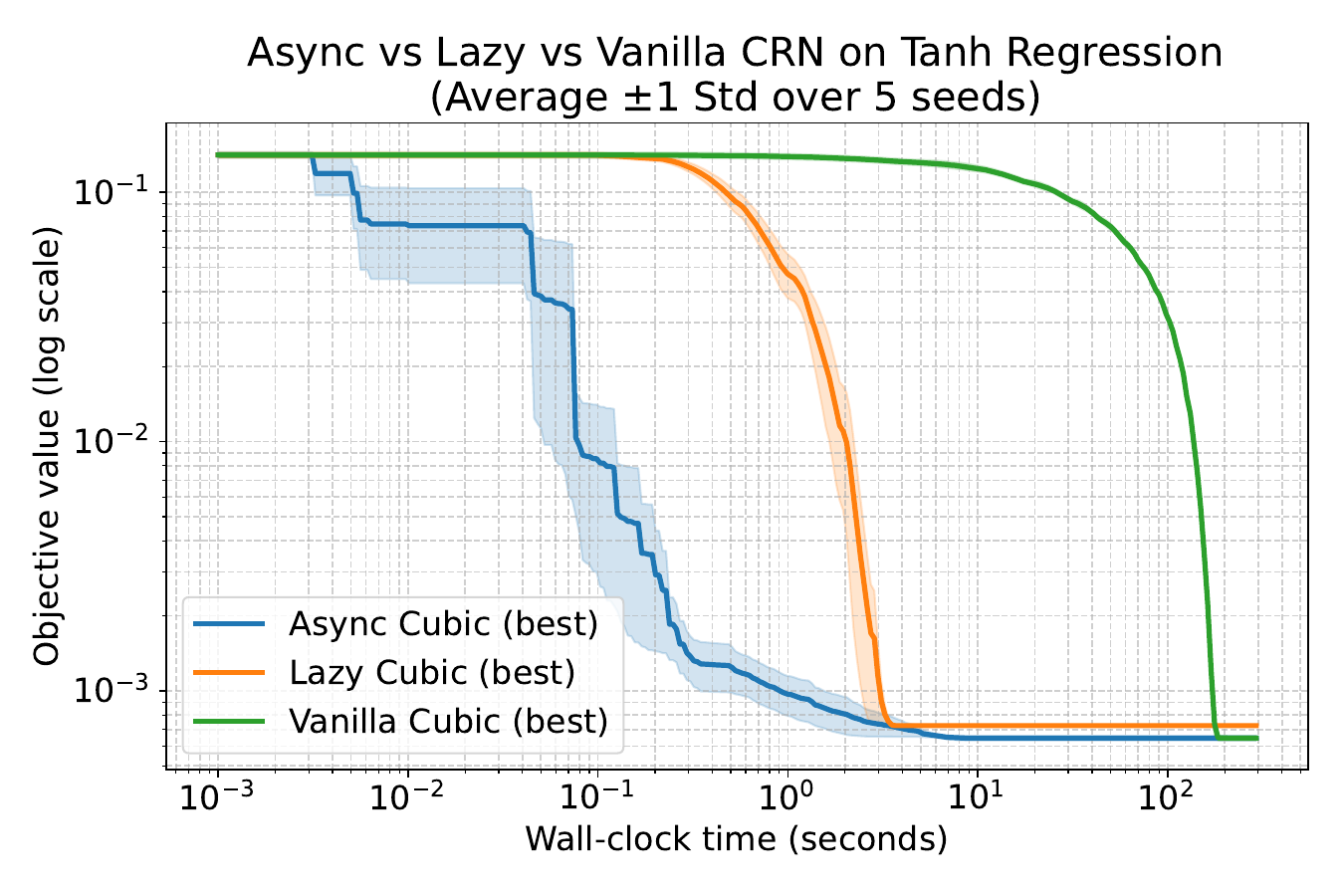}
    \caption{Tanh regression ($d{=}500$): objective versus wall-clock time, averaged over 5 random seeds. The shaded regions denote $\pm 1$ standard deviation.}
    \label{fig:tanh-loss-appendix}
  \end{minipage}
\end{figure}

\paragraph{Hyperparameters used in \cref{fig:tanh-loss-appendix}.}
Async: $\rho = 1.0\times 10^2$. Lazy: $\rho = 1.0 \times 10^1$, $m = 500$. Vanilla: $\rho = 1.0 \times 10^1$. The randomized-SVD oracle introduces additional curvature error.

\subsection{Setting for the LIBSVM experiments}

For the experiments of the main text, we evaluate our methods on the \texttt{a1a} and \texttt{a9a} binary classification datasets from the LIBSVM repository \citep{chang2011libsvm}, with $n = 1605$ and $n = 32561$ training samples. Labels are mapped to $\{0, 1\}$ and features are used without normalization or scaling. The working dimension is the largest feature index present in each training file: $d = 119$ for \texttt{a1a}, whose official schema lists 123 features, and $d = 123$ for \texttt{a9a}. On \texttt{a1a} the Geman--McClure regularization weight is $\lambda = 0.01$. The initial point is drawn as $x_0 \sim \mathcal{N}(0, \sigma^2 I)$ with $\sigma = 0.05$. Consistent with the synthetic experiments, our evaluation protocol proceeds in two steps: we first tune each algorithm, searching $\rho$ over powers of $10$ and then refining around the best value, select the best among the tested configurations, and then we fix these configurations to rerun the experiments independently across 5 random seeds. This allows us to report the average wall-clock convergence and robust standard deviation bands. Because the Hessians in these tasks can be indefinite, we utilize the randomized-SVD low-rank spectral approximation for the cubic subproblem solver.

\paragraph{Geman--McClure regression (\texttt{a1a}).}
We optimize the Geman--McClure objective, which contains a nonconvex penalty, on the \texttt{a1a} dataset ($d=119$). 
\paragraph{Hyperparameters used (\cref{fig:nonconvex-loss}).} 
Async: $\rho = 10$. Lazy: $\rho = 10$, $m = 119$. Vanilla: $\rho = 100$. 

\paragraph{One-layer tanh regression (\texttt{a9a}).}
We optimize the nonconvex one-layer tanh regression objective on the \texttt{a9a} dataset ($d=123$). 
\paragraph{Hyperparameters used (\cref{fig:tanh-loss}).} 
Async: $\rho = 1$. Lazy: $\rho = 1$, $m = 123$. Vanilla: $\rho = 1$.

\paragraph{Dimensionality scaling (Speedup vs.\ $d$).}
To probe the per-step cost ratio $q = \Theta(d + d^2/n)$ (\cref{fig:speedup-vs-d}), we benchmark the steady-state iterations per second (iter/s) across varying problem dimensions. We use the synthetic Geman--McClure regression objective with a fixed sample size of $n=5000$ and sweep the dimension $d \in \{200, 500, 1000, 2000\}$. Each algorithm is executed for a fixed 5-second wall-clock budget. To ensure robustness, the iter/s ratios are computed as an average across 5 random seeds, capturing stable steady-state performance while accounting for system and threading variance.
\subsection{Empirical validation of (A3) cost stability}
\label{app:cost_stability}

\cref{lem:overlap} uses (A3): $\Delta_{i+1} \le \gamma\Delta_i$, where $\Delta_i = a_{i+1} - a_i$ counts gradient iterations between successive dispatches. We measured instead the wall-clock duration of each curvature job, normalized by the iteration cost; across the measured runs, adjacent ratios of this physical quantity stayed within $[0.9, 1.1]$. The two measurements need not agree after rounding to integer counts, or when iteration durations vary, so this physical-time diagnostic does not certify (A3) in the units of the proof, and we do not derive the constants of the theorems from it. It is consistent with predictable dense linear-algebra timing on a non-shared CPU at fixed $d$, and concerns the measured runs only.

\subsection{Measured delays, curvature-oracle error and cost profiling}
\label{app:assumption_diagnostics}

We log the realized delay sequence $\{\tau_t\}$ and Hessian inexactness $\{\|\Ht_t - \hess(x_{t-\tau_t})\|_F\}$ from the Geman--McClure run of \cref{sec:experiments}, plus peak memory across a dimension sweep.

\begin{figure}[h]
  \centering
  \begin{minipage}[b]{.32\linewidth}
    \centering
    \includegraphics[width=\linewidth]{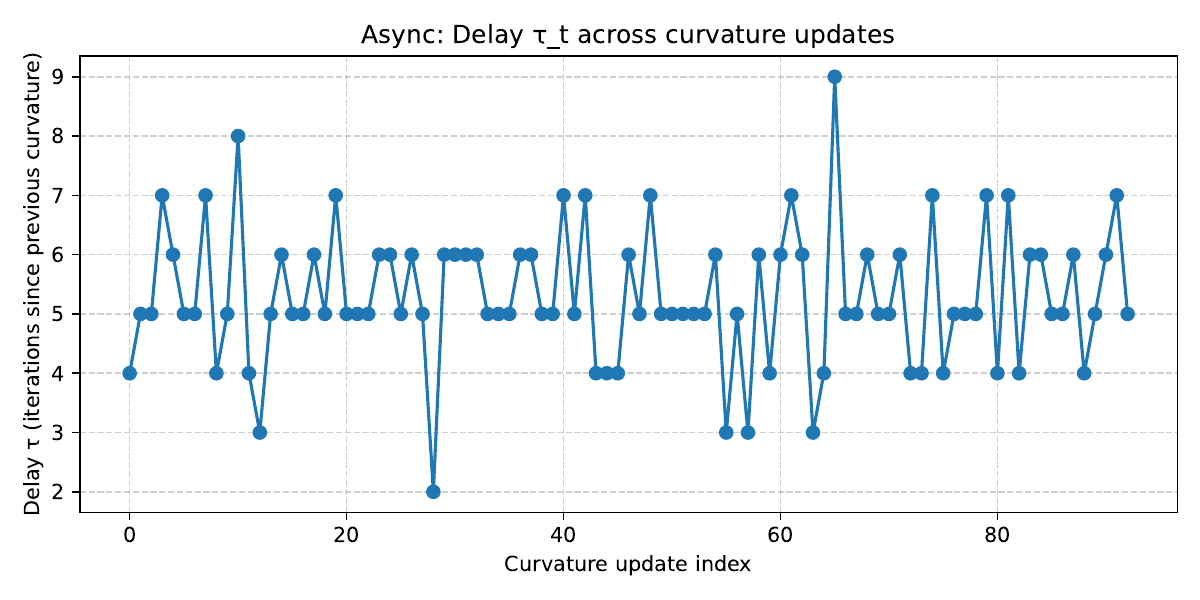}
    \caption{Observed ages at curvature-update times, in $\{2, \dots, 9\}$ over $\sim 90$ updates, with mean $\approx 5.4$. These are not the full per-iteration age sequence $\hat\tau_k$.}
    \label{fig:async-delay}
  \end{minipage}\hfill
  \begin{minipage}[b]{.32\linewidth}
    \centering
    \includegraphics[width=\linewidth]{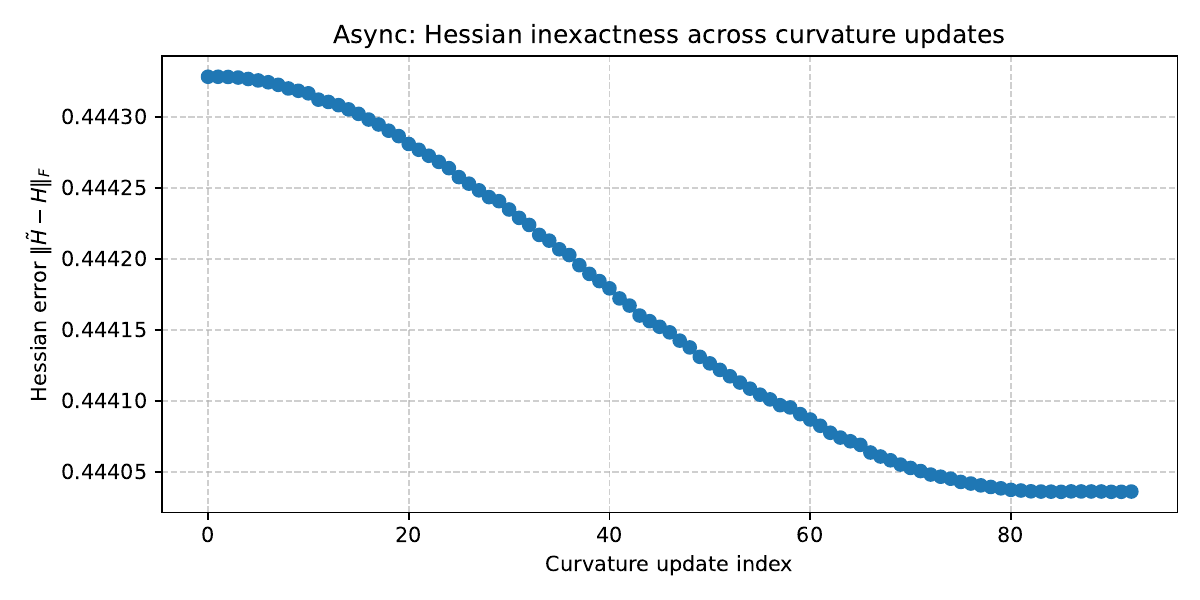}
    \caption{Frobenius-norm error of the curvature oracle on the measured run, in $[0.44403, 0.44430]$. This concerns the later oracle error, relevant to the uniform-error setting of \cref{thm:robust}, not (A5); the narrow spread is not a small relative approximation error.}
    \label{fig:async-inexactness}
  \end{minipage}\hfill
  \begin{minipage}[b]{.32\linewidth}
    \centering
    \includegraphics[width=\linewidth]{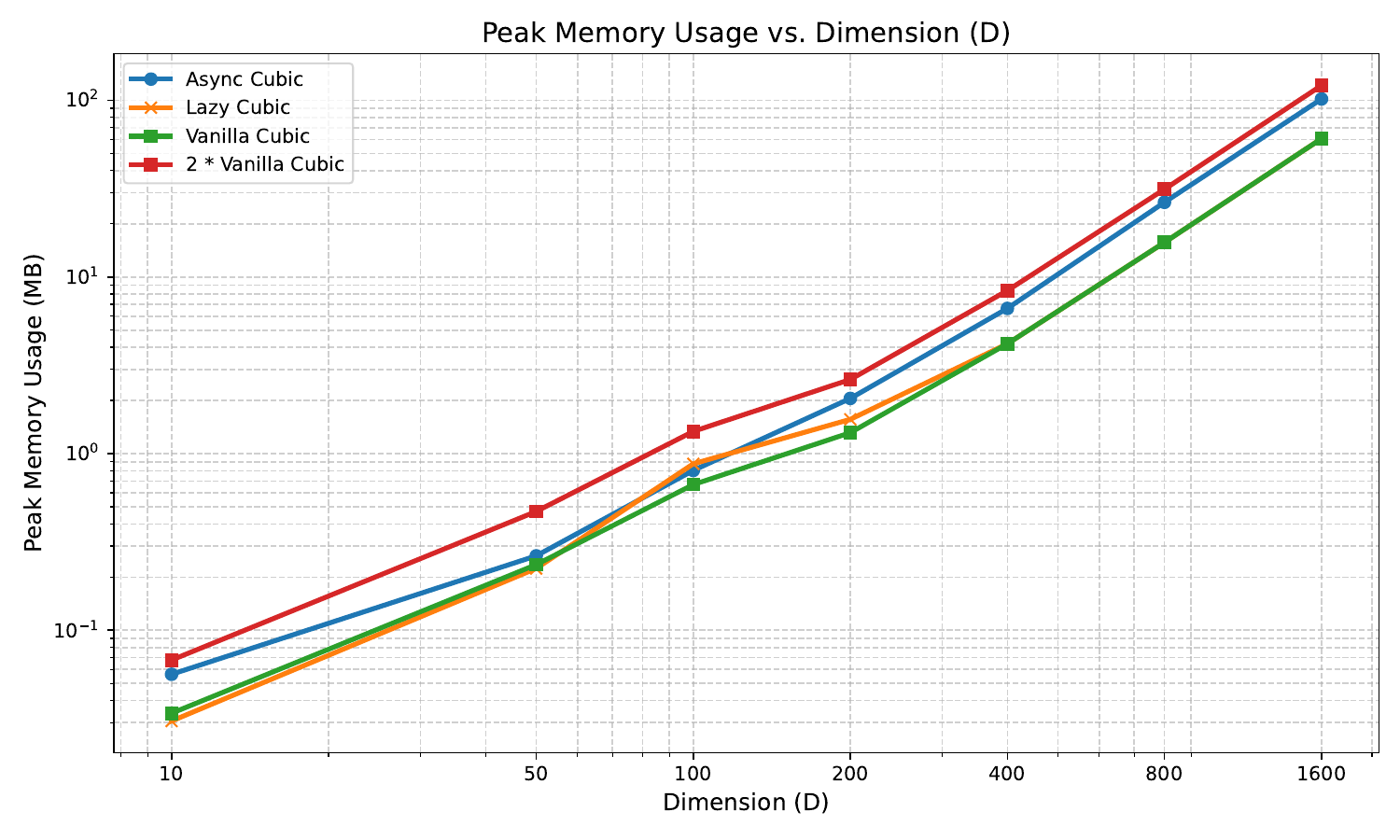}
    \caption{Peak memory of Async tracks $\sim 2\times$ Vanilla (double-buffering of factorizations); $\sim 100$MB at $d=1600$.}
    \label{fig:async-memory}
  \end{minipage}
\end{figure}

\cref{fig:async-delay} reports the ages observed at curvature-update times: they fluctuate over $\{2, 3, \ldots, 9\}$ with no upward trend on this run, which does not by itself establish (A4) for the full per-iteration sequence. On this run, the curvature-oracle error remains in the interval reported in \cref{fig:async-inexactness}; this diagnostic does not isolate the effect of asynchrony on approximation quality, which would require a corresponding synchronous oracle experiment. \cref{fig:async-memory} reports peak memory across $d \in \{10, 50, \ldots, 1600\}$: Async tracks $\sim 2\times$ Vanilla as predicted by the double-buffering design, a deliberate space-time trade-off, modest in absolute terms.

\begin{figure}[t]
  \centering  \includegraphics[width=0.4\linewidth]{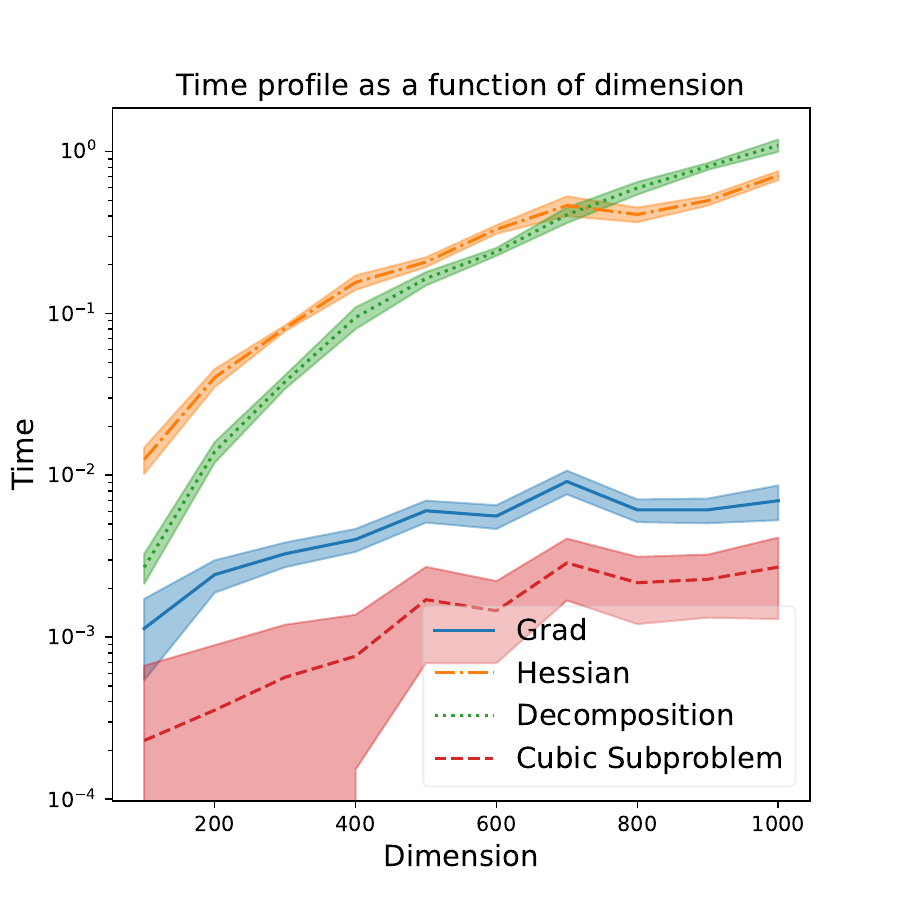}
  \caption{\textbf{Cost profiling:} per-call time in seconds of each component versus $d$ (mean $\pm$ std over seeds), with the Hessian assembled in full at every call. Assembly and factorization dominate the gradient and the cubic solve.}
  \label{fig:time-profile-wide}
\end{figure}

\paragraph{Cost Profiling.}
\Cref{fig:time-profile-wide} profiles, per call, a gradient evaluation, Hessian assembly, a dense Cholesky factorization and a cubic-subproblem solve, as functions of $d$. It is a cost microbenchmark, separate from the randomized-SVD oracle of the convergence runs. Hessian assembly and factorization together dominate the per-iteration cost and grow much faster than the gradient; assembly is the larger of the two up to $d \approx 700$, and factorization beyond.

\end{document}